\documentclass{amsart}
\usepackage{amssymb}
\usepackage{amsmath}
\usepackage{amsthm}
\usepackage{amsfonts}
\usepackage{graphicx}
\usepackage{enumerate}
\usepackage[ansinew]{inputenc}
\usepackage{hyperref}
\usepackage{multirow}
\usepackage{color}
\usepackage{colortbl}   

\usepackage{wrapfig}
\usepackage{floatflt,epsfig} 
\usepackage{fancyhdr}
\usepackage{array}
\usepackage[right]{eurosym}

\newtheorem{theorem}{Theorem}
\newtheorem{def.}[theorem]{Definition}
\newtheorem{prop.}[theorem]{Proposition}
\newtheorem{lem.}[theorem]{Lemma}
\newtheorem{cor.}[theorem]{Corollary}
\newtheorem{conj.}[theorem]{Conjecture}
\newtheorem{example}[theorem]{Example}

\newtheorem{rem.}[theorem]{Remark}

\def\Hil{\mathcal{H}}

\def\be{\begin{equation}}
\def\ee{\end{equation}}
\def\ben{\begin{eqnarray}}
\def\een{\end{eqnarray}}

\def\AAa{{\mathcal A}}

\def\HHmp{{{\mathcal H}_w^p}}

\definecolor{darkviolet}{rgb}{0.58,0,0.83} 

\usepackage{tikz}
\usepackage{enumitem}
\usepackage{comment}
\allowdisplaybreaks[1]
\newcommand{\B}[0]{{\mathcal{B}}}

\newcommand{\R}[0]{{\mathbb{R}}}

\newcommand{\G}[0]{{\mathcal{G}}}
\newcommand{\N}[0]{{\mathbb{N}}}
\newcommand{\Z}[0]{{\mathbb{Z}}}
\newcommand{\A}[0]{{\mathcal{A}}}
\newcommand{\J}[0]{{\mathcal{J}}}
\newcommand{\MS}[0]{{\mathcal{S}}}

\renewcommand{\H}{{\mathcal{H}}}

\newtheorem{proposition}[theorem]{Proposition}
\definecolor{darkviolet}{rgb}{0.58,0,0.83} 

\def\lukas{\color{blue}}

\title{Localized frames without inequalities}
\author{P. Balazs}
\address{(P. B.)  Acoustics Research Institute, Austrian Academy of Sciences,
Dominikanerbastei 16, 1010 Vienna}
\address{Austria and Acoustics, Analysis and AI, Interdisciplinary Transformation University Austria, Science Park 4, 4040 Linz, Austria}
\email{peter.balazs@oeaw.ac.at}
\author{L. K\"ohldorfer}
\address{(L. K.) Acoustics Research Institute, Austrian Academy of Sciences, Dominikanerbastei 16, 1010 Vienna, Austria}
\email{koehli@gmx.at}
\author{M. Speckbacher}
\address{(M. S.) Acoustics Research Institute, Austrian Academy of Sciences, Dominikanerbastei 16, 1010 Vienna, Austria }
\email{michael.speckbacher@oeaw.ac.at}
\date{}

\begin{document}

\begin{abstract}
We consider countable families of vectors in a separable Hilbert space, which are mutually localized with respect to a fixed localized Riesz basis. We prove the equivalence of the frame property and nine conditions that do not involve any inequalities. This is done by studying the properties of their frame-related operators on the co-orbit spaces generated by the reference Riesz basis. We apply our main result to the setting of shift-invariant spaces and obtain new conditions for stable sets of sampling. 
\end{abstract}

\maketitle

\section{Introduction}

A frame is a natural generalization of an orthonormal basis for a Hilbert space that allows for redundant representations. This makes frames particularly well suited to applications such as signal processing.
A countable family $(\psi_k)_{k\in X}$ of vectors in a separable Hilbert space $\Hil$ is called a frame, if there exist $0< A \leq B$ such that \begin{equation}\label{eq:frame}
A \, \| f \|^2 \leq
  \sum_{k\in X} \vert \langle f  ,  \psi_k \rangle \vert^2 \leq
  B \, \| f \|^2, \qquad  f\in\Hil.
\end{equation} 
 On the Hilbert space level, the frame inequalities (\ref{eq:frame}) admit several well-known reformulations in terms of the canonical frame-related operators, such as the analysis operator being bounded, injective, and having closed range, or the synthesis operator being bounded and surjective, or the frame operator being boundedly invertible (see \cite{xxlstoeant11} or Section~\ref{Basic Frame Theory}). While these characterizations are fundamental, in concrete situations they often merely recast the frame inequalities in operator language, without making an underlying mechanism or obstruction more transparent. This motivates the question whether the frame property of a given family of vectors can be determined by conditions that do not merely reformulate the frame inequalities.

A fundamental result in this direction is due to Gr\"ochenig \cite{gro07ineq}, who gave fourteen equivalent conditions \emph{without inequalities} for a regular Gabor system with non-zero window in the modulation space $M^1(\R^d)$ being a frame. 
The main result of \cite{gro07ineq} is that any such Gabor system $\G = \G(g,\Lambda)$ is a frame for $L^2(\mathbb{R}^d)$ if and only if the associated frame-related operators satisfy certain mapping properties on modulation spaces, for example, if the analysis operator $C_{\G}: M^\infty(\mathbb{R}^{d}) \to \ell^{\infty}(\Lambda)$ is injective. The proof of the latter relies on deep results 
of Gabor systems over lattices, such as the Ron-Shen duality principle \cite{ron-shen}, Wiener's Lemma for twisted convolutions \cite{GrLe04}, and a result on linear independence of time-frequency shifts 
\cite{gro07ineq}.

In many settings, however, frame theoretic results are ultimately not the consequence of an underlying group structure, but follow from localization properties of the frame, see, e.g., kernel theorems (compare \cite{xxlgrospeck19,feichtinger80cras} with \cite{xxlbysp24}), invertibility results of Gabor frame operators on modulation spaces \cite{groe04,GrLe04,GROCHENIG2015388}, or a partial generalization \cite[Theorem~5.1]{GROCHENIG2015388} of \cite[Theorem~3.1]{gro07ineq}. A frame $(\psi_k)_{k\in X}$ is said to be \emph{localized}, if its Gram matrix has some suitable off-diagonal decay, or, more precisely, if it is contained in a \emph{spectral algebra} (see \cite{forngroech1} or Definition~\ref{def:spectral}). 
To a localized frame $\psi$ one can then associate a family of Banach spaces, the so-called \emph{co-orbit spaces} $\H^p(\psi)\ (1\leq p\leq \infty)$ which are defined in terms of $\ell^p$-membership of the analysis coefficients $(\langle f,\psi_k \rangle)_{k\in X}$, see Section~\ref{sec:localized-frames} for further details.

Returning to Gabor analysis, we note that the property $g\in M^1(\mathbb{R}^d)$ of the window function can be seen as a localization assumption. Indeed, in the case $\Lambda = \mathbb{Z}^{2d}$, we note that $\G(g,\mathbb{Z}^{2d})$ is localized with respect to the unweighted Baskakov-Gohberg-Sj\"ostrand algebra $\mathcal{C}$ (see Example~\ref{spectralexamples} (3)), and in the case of more irregular index sets $\Lambda$, similar (more general) properties have been shown \cite{GROCHENIG2015388}. Furthermore, if $\G$ is a frame, the modulation spaces $M^p(\R^d)$ coincide with the co-orbit spaces $\Hil^p(\G)$ ($1\leq p \leq \infty$) generated from $\G$ in this case \cite[Theorem~3.2]{GROCHENIG2015388}. Naturally, 
these observations lead to the question, whether abstract localized frames can be characterized in a similar fashion as in \cite{gro07ineq}, namely via properties of their associated frame-related operators on their respective co-orbit spaces. 

In this article we give several equivalent conditions for a suitably localized family \(\psi=(\psi_k)_{k\in X}\) of vectors in an abstract separable Hilbert space \(\mathcal H\) being a frame. A central feature of our approach is that the relevant frame-related operators are considered not only on the Hilbert space level, but simultaneously on the associated co-orbit spaces \(\mathcal H^p(\varphi)\). Passing from boundedness below on one of the corresponding sequence spaces \(\ell^p(X)\) to boundedness below on \(\ell^2(X)\) requires a compatibility property of the underlying matrix algebra. We therefore work with uniformly $p$-stable spectral algebras, meaning that, for matrices in the algebra, boundedness below on one space \(\ell^p(X)\) implies boundedness below on every \(\ell^q(X)\), \(1\leq q\leq\infty\). All standard spectral algebras considered in this article satisfy this property, a result which may be of independent interest; see Theorem~\ref{thm:p-stable-examples}.

Let $C_\psi$, $D_\psi$ and $S_\psi$ denote the analysis, synthesis and frame operator respectively, $G_\omega = [\langle \omega_l, \omega_k \rangle ]_{k,l \in X}$ the Gram matrix, 
and set $\text{Ran}_p(C_\psi):=C_\psi(\H^p(\varphi))$ and $\text{Ran}_p(D_\psi):=D_\psi(\ell^p(X))$, see Section~\ref{Basic Frame Theory} for the detailed definitions. The formulation of our main result is as follows:

\begin{theorem}\label{thm:main}
Let $\A$ be a uniformly $p$-stable spectral algebra and $\varphi = (\varphi_k)_{k\in X}$ be an intrinsically $\A$-localized Riesz basis in $\Hil$. Let $\psi = (\psi_k)_{k\in X}$ be a family in $\Hil$, which is mutually $\A$-localized with respect to $\varphi$, and let $\omega = (\omega_k)_{k\in X}$ be given by 
$$\omega_k = \sum_{l\in X}\langle \psi_l , \varphi_k \rangle S_{\varphi}^{- {1}/{2}}\varphi_l, \qquad k\in X.$$
Then the following are equivalent.
\begin{enumerate}[label=(\arabic*)]
    \item\label{list:1} $\psi$ is a frame for $\Hil$.
    \item\label{list:2} $S_{\psi}$ is invertible on $\Hil^1(\varphi)$.
    \item\label{list:3} $S_{\psi}$ is invertible on $\Hil^{\infty}(\varphi)$.
          \item\label{list:4}  $C_{\psi}:\Hil^{\infty}(\varphi) \rightarrow \ell^{\infty}(X)$ is injective and  $\emph{Ran}_\infty(C_{\psi})$ as well as $\emph{Ran}_\infty(D_\psi)$ are closed.
       \item\label{list:5}  $D_{\psi}:\ell^1(X) \rightarrow \Hil^1(\varphi)$ is surjective and $\emph{Ran}_1(C_\psi)$ is closed.
         \item\label{list:6}  $D_{\omega}:\ell^{\infty}(X) \rightarrow \Hil^{\infty}(\varphi)$ is injective and $\emph{Ran}_\infty(C_\omega)$ as well as $\emph{Ran}_\infty(D_\omega)$ are closed.
      \item\label{list:7}  $C_{\omega}:\Hil^1(\varphi) \rightarrow \ell^1(X)$ is surjective and $\emph{Ran}_1(D_\omega)$ is closed.
    \item\label{list:8}  $G_{\omega}$ is invertible on $\ell^1(X)$. 
    \item\label{list:9}  $G_{\omega}$ is invertible on $\ell^{\infty}(X)$.
    \item\label{list:10}  $\omega$ is a Riesz sequence in $\H$.
\end{enumerate}
In particular, if any of the above conditions is satisfied, then $\Hil^p(\psi) = \Hil^p(\varphi)$ with equivalent norms ($1\leq p \leq \infty$).
\end{theorem}

At this point, we explicitly emphasize that Theorem \ref{thm:main} complements, but does \emph{not} include Gr\"ochenig's main result in \cite[Theorem~3.1]{gro07ineq} as a special case. We refer to the bullet points below and to Section~\ref{A Comparison to "Gabor Frames Without Inequalities"} for a comparison of our main result to the one in \cite{gro07ineq}. Nevertheless, \cite[Theorem~3.1]{gro07ineq} served as the main inspiration for the present work, and certain parts of the underlying proof strategy (e.g. modulation spaces, duality relations, or the closed range theorem) could indeed be transferred to the abstract setting.   

In what follows, we outline some of our  initial observations that led to our approach and, in particular, the first key differences to the Gabor setting \cite{gro07ineq}, that were decisive for the formulation of our main result:

\begin{itemize}
\item \textbf{Necessity of a localized reference system:} While the modulation spaces $M^p(\R^d)$ can be independently defined via the \emph{short-time Fourier transform} associated with the omnipresent Gaussian (or any other Schwartz function) \cite{gr01}, the definition and many relevant properties of the co-orbit spaces $\Hil^p(\varphi)$ depend on the frame property of some explicitly given localized frame $\varphi$. This is the reason why we relate the family $\psi$ to an intrinsically localized reference frame $\varphi$ generating the co-orbit spaces $\Hil^p(\varphi)$, instead of working with the co-orbit spaces $\Hil^p(\psi)$ directly.
\item \textbf{The R-dual:} To our knowledge, there is currently no canonical generalization of the Ron-Shen duality principle to abstract frames in the literature. 
The Ron-Shen duality principle connects any Gabor frame on a lattice with its (dual) Gabor Riesz sequence generated by the same window on the adjoint lattice. In Theorem~\ref{thm:main} we consider the family $\omega$, which, in the language of \cite{r-duals}, is a so-called \emph{R-dual of type IV} of $\psi$, and serves as a substitute of the aforementioned Ron-Shen dual. In general, however, an R-dual associated with a Gabor frame on a lattice does not coincide with its Ron-Shen dual, see \cite{Casazza2004,r-duals} for a discussion of this matter. 
\item \textbf{Riesz basis property of the reference frame:} In order to guarantee that the aforementioned R-dual $\omega$ is a Riesz sequence if and only if $\psi$ is a frame, we additionally have to assume that the reference system $\varphi$ is in fact a Riesz basis and not only a frame. Thus, our localization assumption $\psi \sim_\A \varphi$ corresponds to the original localization concept from \cite{groe04} in the case $\A = \J_s$ (see Example~\ref{spectralexamples} (1)), which was followed by several other publications on this topic \cite{ALDROUBI20081667,bacahela06a,Balan2006,xxlgro14,forngroech1}. At the same time, however, it might be very difficult, if not even impossible 
for a Gabor frame to be localized with respect to a localized Riesz basis (see also \cite[Section~5]{groe04} for a discussion of this issue). 
\item \textbf{Closed range conditions:} Finally, several closed range conditions appear in our setting but not in the Gabor setting \cite{gro07ineq}. We show in Example~\ref{counterexample} and Remark~\ref{counterexample2}, that each of the closed range conditions appearing in Theorem~\ref{thm:main} is in fact necessary in our setting, since otherwise, the statement is simply not true. For Gabor frames, however, these conditions are automatically satisfied, see Example~\ref{ex:gabor}. 
\end{itemize} 

\


\noindent The usefulness of Theorem~\ref{thm:main} lies in the possibility of choosing the condition that is best adapted to the structure of a given problem. We illustrate this by considering the following simple example.

Let $(\delta_k)_{k\in \mathbb{N}}$ be the canonical orthonormal basis for $\Hil = \ell^2(\mathbb{N})$ and consider the family $\psi = (\psi_k)_{k\in \mathbb{N}}$ given by 
$$\psi_k := \delta_k+\delta_{k+1}, \qquad k\in\mathbb{N}.$$
As observed in \cite[Example~5.4.6]{ole1n}, $\psi$ is a minimal and complete Bessel sequence, but not a frame. There, the failure of the lower frame inequality is demonstrated by constructing a sequence of unit vectors $(f_n)_{n\in\mathbb N}$ such that $\Vert C_\psi f_n \Vert_{\ell^2}\to0$ as $n\to \infty$. Equivalently, although the analysis operator $C_\psi:\ell^2(\mathbb N)\to\ell^2(\mathbb N)$ is injective by completeness of $\psi$, its range is not closed. Theorem~\ref{thm:main} offers a simple alternative: Taking $\varphi=(\delta_k)_{k\in\mathbb{N}}$ as the reference Riesz basis, we obtain $\Hil^p(\varphi) = \ell^p(\mathbb{N})$; and all localization assumptions are clearly satisfied since $G_{\varphi}$ is the identity matrix and $G_{\psi,\varphi}$ is banded (e.g., take the uniformly $p$-stable algebra \(\mathcal A=\mathcal J_2\); see Theorem~\ref{thm:p-stable-examples}). Since $[C_\psi c]_k=c_k+c_{k+1}$ ($c\in\ell^\infty(\mathbb N)$), we find that the alternating sequence $c=( (-1)^k )_{k\in\mathbb N}$ belongs to the kernel of $C_\psi:\ell^\infty(\mathbb{N})\to\ell^\infty(\mathbb{N})$. Condition \ref{list:4} of Theorem~\ref{thm:main} implies that $\psi$ is not a frame, and the non-trivial kernel element $c\in \ell^{\infty}(\mathbb{N} )\setminus \ell^2(\mathbb{N})$ makes this obstruction transparent.

The above example illustrates one of the advantages of Theorem~\ref{thm:main}: A failure of the lower frame inequality, which on the Hilbert space level can be made visible only through a sequence of unit vectors whose analysis coefficients tend to zero, can be alternatively described by an explicit non-trivial kernel element contained in the larger space $\Hil^{\infty}(\varphi)$. More generally, once the localization assumptions are verified, one may use whichever of the equivalent conditions is easiest to verify in the given problem.

As a more substantial application of our main result, we consider stable sampling in shift-invariant spaces. We show that under certain assumptions on the sampling set and the generator of these shift-invariant spaces, all assumptions of Theorem \ref{thm:main} are satisfied, so that we obtain new necessary and sufficient conditions for stable sets of sampling. We postpone the detailed discussion of this example to Section~\ref{Shift-Invariant Spaces}.

Theorem~\ref{thm:main} relies heavily on the theory of localized frames and their associated co-orbit spaces. Often, in the theory of abstract localized frames, the frame property of a given family in $\Hil$ is assumed, and then conclusions are derived from additional localization assumptions on the given frame. 
The philosophy here is, in a sense, the reverse: We start with a given suitably localized family $(\psi_k)_{k\in X}$ in $\Hil$ and then characterize the frame property by operator-theoretic considerations on certain co-orbit spaces. To the best of our knowledge, this approach is novel in the setting of an abstract Hilbert space. 

The proof of Theorem \ref{thm:main} is carried out within the two scales 
$$\Hil^1(\varphi) \subset \Hil \subset \Hil^{\infty}(\varphi) \quad \qquad \text{and} \qquad \quad \ell^1(X) \subset \ell^2(X) \subset \ell^{\infty}(X),$$
which naturally fit into the framework of abstract Banach Gelfand
triples \cite{feilucor08}, see also \cite{haubakoe25}. It combines the duality relations of the analysis and synthesis operators on these spaces with the closed range theorem. Another key ingredient is the proof that the R-dual $\omega$ inherits the localization properties from $\psi$. We obtain this fact by utilizing the holomorphic functional calculus for the Banach algebras $\A$ and $\B(\ell^2(X))$ together with the spectral invariance of $\A$ in $\B(\ell^2(X))$.  As a consequence, this allows us to consider the Gram matrix $G_{\omega}$ not only as an operator on $\ell^2(X)$, but also on $\ell^1(X)$ and $\ell^{\infty}(X)$. To complete the circular chain of implications, we apply uniform $p$-stability to the cross-Gram matrix \(G_{\widetilde{\varphi},\omega}\). Condition \ref{list:7} implies that this matrix is bounded from below on \(\ell^1(X)\). Its boundedness from below on \(\ell^2(X)\) then follows from uniform $p$-stability and yields that \(\omega\) is a Riesz sequence. This is the only step in the proof that requires uniform $p$-stability; in fact, only the implication from \(\ell^1\)-stability to \(\ell^2\)-stability is used.

\

This article is structured as follows. In Section~\ref{Preliminaries and Notation} we fix our notation and summarize some relevant concepts and results from the literature. In particular, we introduce uniformly $p$-stable spectral algebras and present several concrete examples. The proof of this result, which may also be of independent interest, is given in the appendix. In Section~\ref{Auxiliary Results} we prove some auxiliary results, which we utilize in Section~\ref{Proof of main}, where we prove our main result. In Section~\ref{A Comparison to "Gabor Frames Without Inequalities"} we compare Theorem~\ref{thm:main} to \cite[Theorem~3.1]{gro07ineq} and give concrete examples demonstrating the necessity of the closed range conditions. Finally, Section~\ref{Shift-Invariant Spaces} is devoted to applying our main result in the setting of shift-invariant spaces.

\section{Preliminaries and Notation}\label{Preliminaries and Notation}

Throughout this manuscript $\Hil$ denotes a separable Hilbert space and $X$ a countable (index) set. If not explicitly stated otherwise, $\Vert \cdot \Vert$ and the bracket $\langle \cdot,\cdot\rangle$ (without index) denote the norm and the inner product on $\Hil$, respectively. Whenever we talk about the dual space $B'$ of a Banach space $B$, we refer to the \emph{anti-linear} (topological) dual space of $B$. In particular, the identification $(\ell^p(X))' = \ell^q(X)$ (for $1\leq p < \infty$, $\frac{1}{p}+\frac{1}{q} = 1$) is based on the sesquilinear form 
$$\big\langle (c_k)_{k\in X}, (d_k)_{k\in X}\big\rangle_{\ell^p\times \ell^q} = \sum_{k\in X}c_k \overline{d_k}$$
in order to be compatible with the Hilbert space duality on $\ell^2(X)$.

\subsection{Basic Frame Theory}\label{Basic Frame Theory}

A countable family $\varphi : = (\varphi_k )_{k \in X}$ of vectors in 
$\H$ is called a \emph{frame}, 
if there exist constants $A,B>0$ such that 
\begin{equation}\label{framedef}
A \, \| f \|^2 \leq
  \sum_{k\in X} \vert \langle f  ,  \varphi_k \rangle \vert^2 \leq
  B \, \| f \|^2, \qquad  f\in\Hil.
\end{equation} 
  For a comprehensive introduction to frames in Hilbert space we refer the reader to \cite{ole1n}.
If the upper, but not necessarily the lower inequality in (\ref{framedef}) is satisfied for some prescribed $B>0$ and all $f\in \Hil$, then $\varphi$ is called a \emph{Bessel sequence}. Whenever $\varphi$ is a Bessel sequence, both the \emph{analysis operator}, defined as
\begin{equation*}
C_{\varphi} : \H \rightarrow  \ell^2(X) ,\qquad C_\varphi f= \big(\langle f   ,   \varphi_k \rangle\big)_{k \in X},
\label{eq:Analysenoperator}
\end{equation*}  
and the \emph{synthesis operator}, defined as
\begin{equation*}
D_\varphi : \ell^2(X) \rightarrow \H ,\qquad D_\varphi c =  \sum_{k \in X} c_k \varphi_k,
\label{eq:Synthesenoperator}
\end{equation*}
are bounded and adjoint to one another, where the latter series converges unconditionally in $\Hil$. A combination of the analysis and the synthesis operator yields the \emph{frame operator}
\begin{equation*}
S_\varphi  = D_\varphi C_\varphi: \H \rightarrow \H ,\qquad S_\varphi f = \sum_{k \in X} \langle f   ,   \varphi_k \rangle  \varphi_k.
\label{eq:Rahmenoperator}
\end{equation*}
If $\varphi$ is a frame, then $S_\varphi$ is bounded, positive, self-adjoint and invertible and the family $\widetilde{\varphi}=\big(\widetilde{\varphi}_k \big)_{k \in X} :=\big( S_{\varphi}^{-1} \varphi_k \big)_{k \in X}$ is another frame, called  the \emph{canonical dual frame}, that satisfies 
\begin{equation}\label{framerec}
f=D_\varphi C_{\widetilde{\varphi}}f=\sum_{k \in X} \langle f   ,   \widetilde{\varphi}_k \rangle  \varphi_k=D_{\widetilde{\varphi}}C_\varphi f=\sum_{k \in X} \langle f   ,   \varphi_k \rangle  \widetilde{\varphi}_k, \qquad  f\in \Hil.
\end{equation}
More generally, if $\varphi^d = (\varphi^d)_{k\in X}$ is another family in $\Hil$ such that $D_\varphi C_{\varphi^d} = D_{\varphi^d} C_{\varphi} = \mathcal{I}_{\Hil}$, then $\varphi^d$ is called a \emph{dual frame} of $\varphi$. If $\psi,\varphi$ are two countable families of vectors in $\H$ indexed by $X$, then we define the \emph{cross Gram matrix} as
\begin{equation}\label{eq:Kreuz-Gramsche_Matrix}
G_{\psi,\varphi}  = \big[ \langle \varphi_{l} , \psi_k \rangle \big]_{k,l \in X}.
\end{equation}
If $\psi=\varphi$ we write $G_\varphi=G_{\varphi,\varphi}$. In general, a (cross) Gram matrix as in (\ref{eq:Kreuz-Gramsche_Matrix}) might not be bounded on $\ell^2(X)$, but in case $\psi$ and $\varphi$ both are Bessel sequences, $G_{\psi,\varphi}$ defines a bounded operator on $\ell^2(X)$ and 
$$G_{\psi,\varphi}=C_\psi D_\varphi $$
and $G_{\psi,\varphi}^\ast=G_{\varphi,\psi}$. 
In fact, a countable family $\varphi$ in $\Hil$ is a Bessel sequence if and only if its associated Gram matrix $G_{\varphi}$ defines a bounded operator on $\ell^2(X)$ \cite{xxlstoeant11}. A family $\varphi : = (\varphi_k )_{k \in X}\subset \H$   is called a \emph{Riesz sequence} for $\H$, if there exist constants $A,B>0$ such that
$$A\|c\|_{\ell^2(X)}^2\leq \left\Vert \sum_{k\in X}c_k\varphi_k\right\Vert^2\leq B\|c\|_{\ell^2(X)}^2, \qquad  c\in \ell^2(X).$$
Equivalently, $\varphi$ is a Riesz sequence if and only if its associated Gram matrix $G_{\varphi}$ defines a bounded and invertible operator on $\ell^2(X)$ \cite{xxlstoeant11}. A Riesz sequence $\varphi$ is called a \emph{Riesz basis}, if it is also \emph{complete}, that is, $\overline{\text{span}}(\varphi_k)_{k\in X} = \H$. In particular, every Riesz sequence is a Riesz basis for $\overline{\text{span}}(\varphi_k)_{k\in X}$. We also note that a Riesz basis is a frame for which both $G_{\varphi}\in \B(\ell^2(X))$ and $S_{\varphi}\in \B(\Hil)$ are invertible. For any Riesz basis $\varphi$ the family $(S_\varphi^{-1/2}\varphi_k)_{k\in X}$ constitutes an orthonormal basis. 

\subsection{The R-dual sequence}

We will consider the following notion of duality for abstract sequences in $\H$. Let $(e_k)_{k\in X}$ and $(h_k)_{k\in X}$ be two orthonormal bases for $\H$, and let $(\phi_k)_{k\in X}$ be a family of vectors such that $\sum_{n\in X} |\langle \phi_n,e_k\rangle |^2<\infty$ for every $k\in X$.
Then 
\begin{equation}\label{Rdual0}
\omega_k=\sum_{n\in X}\langle \phi_n,e_k\rangle h_n
\end{equation}
converges in $\Hil$ for each $k\in X$ and the family $\omega = (\omega_k)_{k\in X}$ is called an 
\emph{R-dual} or \emph{Riesz-dual sequence} for $\phi$. It was shown in \cite[Propositions 10,11]{Casazza2004} that $\phi$ is a frame for $\H$ if and only if $\omega$ is a Riesz sequence in $\Hil$.

The R-dual $\omega$ that we use in Theorem \ref{thm:main} fits into the framework of (\ref{Rdual0}), since
$$\omega_k = \sum_{l\in X}\langle \psi_l , \varphi_k \rangle S_{\varphi}^{- {1}/{2}}\varphi_l = \sum_{l\in X}\langle S_{\varphi}^{ {1}/{2}}\psi_l , S_{\varphi}^{- {1}/{2}}\varphi_k \rangle S_{\varphi}^{- {1}/{2}}\varphi_l$$
and $(S_{\varphi}^{- {1}/{2}}\varphi_k)_{k\in X}$ is an orthonormal basis. Note that $(\psi_k)_{k\in X}$ is a frame if and only if $(S_{\varphi}^{ {1}/{2}}\psi_k)_{k\in X}$ is a frame, because $S_{\varphi}^{ {1}/{2}} \in \B(\Hil)$ is boundedly invertible.

\subsection{Localized Frames}\label{sec:localized-frames}

As mentioned in the previous section, whenever $\varphi$ is a Bessel sequence, its associated Gram matrix $G_{\varphi} = [\langle \varphi_l,\varphi_k\rangle]_{k,l\in X}$ defines a bounded operator on $\ell^2(X)$. Boundedness on $\ell^2(X)$ is equivalent to certain conditions on the matrix entries $\langle \varphi_l,\varphi_k \rangle$, as described, e.g., in \cite[Theorem 4.16]{Maddox:101881}. The idea of localization is to impose (much stronger) off-diagonal decay conditions on a given (cross) Gram matrix. Off-diagonal decay of matrices is often described by Banach matrix algebras, which motivated the authors in \cite{forngroech1} to consider the following type of matrix algebras.

\begin{def.}\label{def:spectral}
A unital Banach *-algebra $\mathcal{A}$ of $\mathbb{C}$-valued matrices $A=[A_{k,l}]_{k,l\in X}$ is called a \emph{solid spectral matrix algebra} \cite[Section~2.3]{forngroech1}, or simply \emph{spectral algebra}, if
\begin{enumerate}
\item [(A1)] $\A \subset \bigcap_{1\leq p\leq \infty} \B(\ell^p(X))$, i.e., each $A \in \mathcal{A}$ defines a bounded operator on $\ell^p (X)$ for each $1\leq p\leq \infty$.
\item[(A2)] $\mathcal{A}$ is \emph{inverse-closed} in $\mathcal{B}(\ell^2 (X))$, i.e., if $A \in \mathcal{A}$ is invertible in $\mathcal{B}(\ell^2(X))$, then its bounded inverse $A^{-1} \in \B(\ell^2(X))$ is contained in $\mathcal{A}$ as well.
\item[(A3)] $\mathcal{A}$ is \emph{solid}, i.e., if $A = [A_{k,l}]_{k,l\in X} \in \mathcal{A}$ and $\vert B_{k,l} \vert \leq \vert A_{k,l} \vert$ for all $k,l \in X$, then $B = [B_{k,l}]_{k,l\in X} \in \mathcal{A}$ and $\Vert B \Vert_{\mathcal{A}} \leq \Vert A \Vert_{\mathcal{A}}$. 
\end{enumerate}
\label{eq:Spektrale-Matrix-Algebra}
\end{def.}

\begin{rem.}\label{holfunccalc}
Condition (A2) is equivalent to \emph{spectral invariance} of $\A$ in $\B(\ell^2(X))$, that is, $\sigma_{\A}(A) = \sigma_{\B(\ell^2(X))}(A)$ for all $A\in \A$ \cite[Lemma 2.4]{gr10-2}, where $\sigma_{\A}(A)$ and $\sigma_{\B(\ell^2(X))}(A)$ denote the spectrum of $A$ in the Banach algebras $\A$ and $\B(\ell^2(X))$ respectively. As a consequence, the holomorphic (Riesz) functional calculi (see e.g. \cite[Chapter~VII]{conw1}) of $\A$ and $\B(\ell^2(X))$ coincide. More precisely, if 
$A\in \A$, $U\subseteq \mathbb{C}$ is an open neighborhood of (a connected component of) $\sigma_{\A}(A) = \sigma_{\B(\ell^2(X))}(A)$, $\Gamma \subseteq U$ is a contour of $\sigma_{\A}(A) = \sigma_{\B(\ell^2(X))}(A) $, and if $z\mapsto f(z)$ is holomorphic on $U$, then 
$$f(A) = \frac{1}{2\pi i}\int_{\Gamma} f(z)(z - A)^{-1} \, dz$$
defines the same element in $\A$ and $\B(\ell^2(X))$. 
\end{rem.}

Strictly speaking, the authors in \cite{forngroech1} considered algebras satisfying properties (A2)$-$(A3) plus the slighly less restrictive assumption (A$1^\ast$), where
\begin{itemize}
 \item [(A$1^\ast$)] $\A \subset \B(\ell^2(X))$, i.e., each $A \in \mathcal{A}$ defines a bounded operator on $\ell^2 (X)$. 
\end{itemize} Here we assume for simplicity, that a spectral algebra  satisfies (A1), since all of the following classical examples of spectral algebras do so anyway, see below. In the language of \cite[Definition~2]{forngroech1}, property (A1) can be rephrased to saying that the constant weight $1$ is \emph{$\A$-admissible with critical index $p_0=1$}. Also note that in \cite{lifting} the authors coined the term ARI-algebra for algebras satisfying conditions (A1) and (A2).  

\begin{example}\label{spectralexamples}
Let $X\subset \mathbb{R}^d$ be a relatively separated set. Then the following classes of matrices are examples of spectral algebras in the sense of Definition~\ref{def:spectral}, see also \cite{KoeBa25} for their $\B(\Hil)$-valued analogs.

\noindent\emph{(1)} For every $s>d$, the \emph{Jaffard algebra} $\J_s = \J_s(X)$ of matrices $A = [A_{k,l}]_{k,l \in X}$ for which 
$$\Vert A \Vert_{\J_s} := \sup_{k,l\in X} \vert A_{k,l} \vert (1+\vert k-l \vert)^s < \infty$$
is a spectral algebra \cite{ja90}.

\noindent\emph{(2)} For every continuous function $\nu:\mathbb{R}^d\to [0,\infty)$ which
\begin{itemize}
    \item[(2a)] is of the form $\nu(x) = e^{\rho(\Vert x \Vert)}$, where $\rho:[0,\infty) \longrightarrow [0,\infty)$ is a continuous and concave function with $\rho(0) = 0$, and $\Vert   \cdot   \Vert$ any norm on $\mathbb{R}^d$, and
    \item[(2b)] satisfies the GRS-condition (Gelfand-Raikov-Shilov condition \cite{GRS64})
    $$\lim_{n\rightarrow \infty} (\nu(nz))^{\frac{1}{n}} = 1, \qquad z\in \mathbb{R}^{d},$$ and
    \item[(2c)] satisfies the weak growth condition \cite{groelei06,barnes87}
$$\nu(x) \geq C (1+\vert x \vert)^{\delta} \qquad \text{for some } \delta\in (0,1], C>0,$$
\end{itemize}
the weighted \emph{Schur algebra} $\MS^1_{\nu} = \MS^1_{\nu}(X)$ of matrices $A = [A_{k,l}]_{k,l \in X}$ for which
$$\Vert A \Vert_{\MS^1_{\nu}} := \max \left\lbrace \sup_{k\in X} \sum_{l\in X} \vert A_{k,l} \vert \nu(k-l) \, , \, \sup_{l\in X} \sum_{k\in X} \vert A_{k,l} \vert \nu(k-l)\right\rbrace < \infty$$
is a spectral algebra \cite{groelei06,sun}.

\noindent\emph{(3)} For every continuous function $\nu:\mathbb{R}^d\to [0,\infty)$ which
\begin{itemize}
    \item[(3a)] is \emph{submultiplicative}, i.e., $\nu(x+y)\leq \nu(x)\nu(y)$ for every $x,y\in \mathbb{R}^d$, and
    \item[(3b)] is \emph{symmetric}, i.e., $\nu(-x) = \nu(x)$ for every  $x\in \mathbb{R}^d$, and
    \item[(3c)] satisfies the GRS-condition,
\end{itemize}
the weighted \emph{Baskakov-Gohberg-Sj\"ostrand algebra} $\mathcal{C}_{\nu} = \mathcal{C}_{\nu}(\mathbb{Z}^d)$ of matrices $A=[A_{k,l}]_{k,l\in \mathbb{Z}^d}$ for which
$$\Vert A \Vert_{\mathcal{C}_\nu} := \sum_{l\in \mathbb{Z}^d} \sup_{k\in \mathbb{Z}^d} \vert A_{k,k-l} \vert \nu(l) <\infty$$
is a spectral algebra \cite{Baskakov1990,Gohberg89,Sjostrand19941995}.

\noindent\emph{(4)} Anisotropic versions of each of the aforementioned algebras are further examples of spectral algebras \cite[Section~6]{KoeBa25}.
\end{example}  

As already mentioned in the introduction, there is one instance in the proof of Theorem \ref{thm:main}, where an additional uniform $p$-stability assumption of the given spectral algebra is used. 

\begin{def.}\label{def:p-stability}\cite{ShiSun09}
Let \(\mathcal A\) be a set of matrices satisfying (A1) and let \(A\in\mathcal A\). For \(1\leq p\leq\infty\), define the lower stability bound of \(A\) on \(\ell^p(X)\) by
\[
        \lambda_p(A)
        :=
        \inf_{0\neq c\in\ell^p(X)}
        \frac{\|Ac\|_{\ell^p(X)}}{\|c\|_{\ell^p(X)}}.
\]
Then \(\mathcal A\) is called \emph{uniformly $p$-stable} if, for every \(A\in\mathcal A\),
\[
        \lambda_p(A)>0
        \text{ for some }p\in[1,\infty]
        \quad\Longrightarrow\quad
        \lambda_q(A)>0
        \text{ for every }q\in[1,\infty].
\]
\end{def.}

To the best of our knowledge, it is not known whether properties \((\mathrm{A1})\)--\((\mathrm{A3})\) imply uniform $p$-stability. However, all concrete examples of spectral algebras mentioned in Example \ref{spectralexamples} are indeed uniformly $p$-stable, see below. We give a full proof of this result in the appendix.

\begin{theorem}
\label{thm:p-stable-examples}
The spectral algebras $\mathcal{J}_s(X)$, $\MS^1_{\nu}(X)$ and $\mathcal{C}_{\nu}(\mathbb{Z}^d)$, as well as their anisotropic versions, are uniformly $p$-stable.
\end{theorem}

For a scalar matrix $A= [A_{k,l}]_{k,l\in X}$, we define the complex conjugate $\overline{A}$ of $A$ and the absolute value $\vert A\vert$ of $A$ pointwise, i.e., by $\overline{A} := \big[\, \overline{A_{k,l}}\, \big]_{k,l\in X}$ and $\vert A\vert := \big[\, \vert A_{k,l}\vert\, \big]_{k,l\in X}$. We explicitly note that by the solidity condition (A3), it holds 
$$A\in \A \qquad \Longleftrightarrow \qquad \overline{A}\in \A \qquad \Longleftrightarrow \qquad \vert A \vert \in \A .$$
We will use this fact several times without further reference. 

\

Having specified the properties required of the underlying matrix algebra, we now turn to localized families of vectors. 

\begin{def.}\label{eq:Localisation_2_}\cite[Definition~3]{forngroech1}
Let $\A$ be a spectral algebra. Two countable families $\psi = (\psi_k)_{k\in X}$ and $\varphi = (\varphi_k)_{k\in X}$ are called \emph{mutually$\mathcal{A}$-localized}, in short $\psi\sim_\A \varphi$, if $G_{\psi,\varphi}\in\mathcal{A}$. Similarly, $\psi$ is called (intrinsically) $\A$-localized if $\psi\sim_\mathcal{A}\psi$, i.e., $G_{\psi}\in \A$. 
\end{def.}

Note that by involution in $\A$ it holds that $\psi \sim_{\A} \varphi$ if and only if $\varphi \sim_{\A} \psi$. We will use this fact without further reference.

One of the main results in \cite{forngroech1} states that localization of a frame is preserved under canonical duality. More precisely:

\begin{theorem}\label{dualloc}\cite[Theorem~3.6]{forngroech1}
Let $\varphi$ be a frame in $\Hil$. If $\varphi \sim_{\A} \varphi$ then also $\widetilde{\varphi} \sim_{\A} \widetilde{\varphi}$ and $\varphi \sim_{\A} \widetilde{\varphi}$.   
\end{theorem}

The latter theorem immediately implies the following simple lemma, which we will also implicitly apply several times without further notice. For completeness we present the simple proof.

\begin{lem.}\label{easy}
Let $\varphi = (\varphi_k)_{k\in X}$ an $\A$-localized frame for $\Hil$ and $\psi = (\psi_k)_{k\in X}$ be any countable family in $\Hil$ such that $\psi \sim_{\A} \varphi$. Then it holds $\psi \sim_{\A} \widetilde{\varphi}$ and $\psi \sim_{\A} \psi$.
\end{lem.}

\begin{proof}
By assumption, $G_{\psi,\varphi}\in \A$ and by Theorem~\ref{dualloc}, $G_{\widetilde{\varphi}}\in \A$. Thus, by frame reconstruction and the algebra property of $\A$, $G_{\psi,\widetilde{\varphi}} = G_{\psi,\varphi} G_{\widetilde{\varphi}}\in \A$. Similarly, the latter implies that $G_{\psi} = G_{\psi,\varphi} G_{\widetilde{\varphi},\psi} \in \A$.   
\end{proof}

A remarkable feature of a localized frame $\varphi$ in $\Hil$ is that the reconstruction formula (\ref{framerec}) does not only converge in the underlying Hilbert space $\Hil$, but also in a whole chain of associated Banach spaces defined below. 

\begin{def.}\label{Hpdef}\cite[Definition~5]{forngroech1}
Let $\mathcal{A}$ be a spectral algebra, $\varphi$ be a frame for $\Hil$, and $\varphi^d$ be some dual frame of $\varphi$ such that $\varphi \sim_{\A} \varphi^d$. Let
\begin{equation}
\begin{split}
\H^{00}(\varphi) := D_{\varphi}(\ell^{00}(X)) = \left\lbrace  \sum_{k \in X} c_k   \varphi_k :  \    c\in \ell^{00} (X) \right\rbrace
\end{split}
\label{eq:H_00}
\end{equation} 
be the space of finite linear combinations of elements from $\varphi$. For $1\le p<\infty,$ the \emph{co-orbit space} $\H^p (\varphi^d,\varphi)$ is defined as the norm completion of $\H^{00}(\varphi)$ with respect to the norm $\|f\|_{\H^p(\varphi^d,\varphi)}:=\big\| C_{\varphi^d} f \big\|_{\ell^p(X)}$.
\end{def.}

Under the above assumptions, the spaces are well-defined non-trivial Banach spaces with $\H^p (\varphi^d,\varphi) \subseteq \H^q (\varphi^d,\varphi)$ for all $1\leq p\leq q <\infty$ and $\H^2 (\varphi^d,\varphi) = \Hil$ with equivalent norms. Furthermore, the analysis operator $C_{\widetilde\varphi}:\Hil^{00}(\varphi)\to \ell^{00}(X)$ extends via density to an isometry $C_{\widetilde\varphi}:\Hil^{p}(\varphi)\to \ell^{p}(X)$ ($1\leq p < \infty$). 

\begin{rem.}\label{indep}
The definition of $\H^p (\varphi^d,\varphi)$ seems to depend heavily on the dual frame pair $(\varphi,\varphi^d)$. However, if $(\psi,\psi^d)$ is another pair of dual frames satisfying $\psi \sim_{\A} \psi^d$, and if $\psi^d \sim_{\A} \varphi$ and $\varphi^d \sim_{\A} \psi$, then $\Hil^p(\varphi^d,\varphi) = \Hil^p(\psi^d,\psi)$ with equivalent norms \cite[Theorem~3.25]{koeba24}. 
\end{rem.}

If $\varphi$ is an intrinsically $\A$-localized frame, then $\varphi \sim_{\A}\widetilde{\varphi}$ and $\widetilde{\varphi} \sim_{\A}\widetilde{\varphi}$ by Theorem~\ref{dualloc}. In particular, due to the above remark, it holds $\Hil^p(\widetilde{\varphi},\varphi) = \Hil^p(\varphi,\widetilde{\varphi})$ with equivalent norms and we write $\Hil^p(\varphi):= \Hil^p(\widetilde{\varphi},\varphi)$ in this case ($1\leq p < \infty$).

The definition of the space $\Hil^{\infty}(\varphi)$ associated to an $\A$-localized frame $\varphi$ in $\Hil$ is more subtle. One could define it as the anti-linear topological dual space of $\Hil^{1}(\varphi)$ (see Theorem~\ref{duality}), or as the completion of $\Hil$ with respect to the Hausdorff locally convex $\sigma(\Hil,\Hil^{00}(\widetilde\varphi))$-topology as done in \cite{forngroech1,groe04}. In this article, we work with the definition from \cite{xxlgro14} (see also \cite{koeba24}) given below. 

\begin{def.}{\label{ERHw}}\cite{xxlgro14,koeba24}
Let $( f_n )_{n\in\N}$ and $( g_n )_{n\in\N}$ be sequences in $\Hil$ and $\varphi$ be a frame for $\Hil$ with canonical dual frame $\widetilde{\varphi}$. Then we call $ ( f_n )_{n\in\N}$ and $( g_n )_{n\in\N}$ equivalent, in short $( f_n  )_{n\in\N} \sim_{\varphi}  (g_n )_{n\in\N}$, if  
$$\lim_{n\rightarrow \infty}  \langle f_n-g_n, \widetilde{\varphi}_k \rangle = 0 \qquad (\forall k\in X) \qquad \text{and} \qquad \sup_{n\in \mathbb{N}}\Vert C_{\widetilde\varphi}(f_n - g_n) \Vert_{\ell^{\infty}(X)} < \infty .$$
\end{def.}
Clearly the relation $\sim_{\varphi}$ defines an equivalence relation on the space of all sequences in $\Hil$.

\begin{def.}\label{Hwinftydef}\cite{xxlgro14,koeba24}
Let $\A$ be a spectral algebra, and $\varphi$ be an $\A$-localized frame for $\Hil$ with canonical dual frame $\widetilde{\varphi}$. Then $\Hil^{\infty}(\varphi)$ is defined to be the space of all equivalence classes $f=[ ( f_n )_{n\in\N}]_{\sim_{\varphi}}$ of sequences $( f_n  )_{n\in\N}$ in $\Hil$ satisfying 
\begin{itemize}
    \item[(i)] $\lim_{n\rightarrow \infty} \langle f_n , \widetilde{\varphi}_k\rangle =: \langle f , \widetilde{\varphi}_k\rangle$ exists for each $k\in X$, 
    \item[(ii)] $\sup_{n\in \mathbb{N}} \Vert C_{\widetilde{\varphi}}f_n \Vert_{\ell^{\infty}(X)} < \infty .$
\end{itemize}
\end{def.}
It is quickly verified that the defining properties (i) and (ii) do not depend on the choice of the representative sequence $(f_n)_{n\in \mathbb{N}}$. Furthermore, condition (ii) implies that
there is a positive constant $C$ such that 
\begin{equation}\label{Cineq}
\vert \langle f_n, \widetilde{\varphi}_k\rangle \vert \leq C \qquad \text{for each } k\in X \text { and each } n\in \mathbb{N}.
\end{equation}
 This implies that $\vert \langle f , \widetilde{\varphi}_k\rangle \vert \leq C$ for all $k\in X$. In particular, 
\begin{equation}\label{Hpinftynorm}
\Vert f \Vert_{\Hil^{\infty}(\varphi)} := \Vert C_{\widetilde{\varphi}} f \Vert_{\ell^{\infty}(X)} = \sup_{k\in X} \left\vert \lim_{n\rightarrow \infty} \langle f_n , \widetilde{\varphi}_k\rangle \right\vert
\end{equation}
is well-defined and, in fact, defines a norm on $\Hil^{\infty}(\varphi)$ with respect to which $\Hil^{\infty}(\varphi)$ is complete. We explicitly emphasize that here $\langle f , \widetilde{\varphi}_k\rangle$ is not an inner product in $\Hil$ but (by abusing notation) the limit of the sequence $\big( \langle f_n,\widetilde{\varphi}_k \rangle  \big)_{n\in\N}$ (for each $k\in X$). Analogously, $C_{\widetilde{\varphi}} f = (\langle f,\widetilde{\varphi}_k\rangle)_{k\in X}$ is understood as a pointwise limit. In this sense, $C_{\widetilde{\varphi}}:\Hil^{\infty}(\varphi)\to \ell^{\infty}(X)$ defines an isometry.

\begin{prop.}\label{Hinftynormequ}
\cite[Theorem~3.8]{forngroech1}, \cite[Corollary~3.27]{koeba24} Let $\A$ be a spectral algebra, and $\varphi =(\varphi_k)_{k\in X}$ be an $\A$-localized frame for $\Hil$ with canonical dual frame $\widetilde{\varphi}$. Then $\Hil^{\infty}(\varphi)$ is a Banach space with respect to the norm $\Vert C_{\widetilde{\varphi}} \cdot \Vert_{\ell^{\infty}(X)}$, and $\Vert C_{\varphi} \cdot \Vert_{\ell^{\infty}(X)}$ is an equivalent norm.
\end{prop.}

The properties of the operators $D_{\varphi}$ and $S_{\varphi}$ on their respective domains are summarized below. 

\begin{prop.}\label{Donto}\cite[Proposition~2.4]{forngroech1},\cite[Proposition~3.20]{koeba24}
Let $\A$ be a spectral algebra, and $\varphi =(\varphi_k)_{k\in X}$ be an $\A$-localized frame for $\Hil$ with canonical dual frame $\widetilde{\varphi}$. Then for each $1\leq p < \infty$, the synthesis operator 
\begin{flalign}\label{syntheisseries}
D_{\varphi}: \ell^{00}(X) \longrightarrow  \Hil^{00}(\varphi), \notag \\
(c_l)_{l\in X} \mapsto \sum_{l\in X} c_l \varphi_l 
\end{flalign}
continuously extends to a bounded and surjective operator $D_{\varphi}: \ell^{p}(X) \to  \Hil^{p}(\varphi)$ and the associated series  (\ref{syntheisseries}) converges unconditionally in $\Hil^{p}(\varphi)$. 

Similarly, if $F_1 \subseteq F_2\subseteq \dots$ is a nested sequence of finite subsets of $X$ such that $\bigcup_{n=1}^{\infty} F_n = X$, then $f_n := \sum_{l\in F_n} c_l \varphi_l$ is contained in $\Hil$ for each $n\in \mathbb{N}$ and 
\begin{equation}\label{synth}
c=(c_l)_{l\in X} \mapsto D_{\varphi}c := \big[( f_n  )_{n\in\N}\big]_{\sim_{\varphi}} = \left[\left( \sum_{l\in F_n} c_l \varphi_l \right)_{n\in\N}\, \right]_{\sim_{\varphi}}
\end{equation}
defines a bounded and surjective operator $D_{\varphi}:\ell^{\infty}(X)\to \Hil^{\infty}(\varphi)$ and 
$D_{\varphi}c$ is independent of the choice of $\lbrace F_n\rbrace_{n=1}^{\infty}$, i.e., (\ref{synth}) converges unconditionally in the $\sigma(\Hil, \Hil^{00}(\widetilde{\varphi}))$-topology.
\end{prop.}

\begin{theorem}\label{Sinvertible}\cite[Theorem~3.8]{forngroech1},\cite[Corollary~3.27]{koeba24}
Let $\A$ be a spectral algebra, and $\varphi =(\varphi_k)_{k\in X}$ be an $\A$-localized frame for $\Hil$ with canonical dual frame $\widetilde{\varphi}$. Then for each $1\leq p \leq \infty$, the frame operator $S_{\varphi} = D_{\varphi}C_{\varphi}:\Hil^p(\varphi) \to \Hil^p(\varphi)$ is bounded and invertible.
\end{theorem}

An important consequence of the latter theorem is that every intrinsically localized frame $\varphi$ in a Hilbert space automatically extends to a Banach frame for each of its associated co-orbit spaces $\Hil^p(\varphi)$, with coefficient space $\ell^p(X)$ \cite{gr91,groe04}. For related Banach-frame constructions in generalized co-orbit theory, we refer to \cite{dahlkeetal08}.

We conclude this section with the following result on the anti-linear topological dual space of $\Hil^p(\varphi)$. 

\begin{theorem}\label{duality}\cite[Proposition~2]{xxlgro14}
Let $\A$ be a spectral algebra, and $\varphi =(\varphi_k)_{k\in X}$ be an $\A$-localized frame for $\Hil$ with canonical dual frame $\widetilde{\varphi}$. Then
$$(\Hil^p(\varphi))' = \Hil^q(\varphi)$$
with equivalent norms for all $1\leq p<\infty$ with $\frac{1}{p}+\frac{1}{q}=1$, where the duality relation is given by 
$$\langle f,g \rangle_{\Hil^p\times \Hil^q} = \langle C_{\widetilde{\varphi}}f,C_{\varphi}g \rangle_{\ell^p\times \ell^q} = \sum_{k\in X} \left[C_{\widetilde{\varphi}}f\right]_k \overline{\left[C_{\varphi}g \right]_k}.$$ 
\end{theorem}

\section{Auxiliary Results}\label{Auxiliary Results}

In this section we prepare the tools that will be important for the proof of our main theorem.

\subsection{Frame-related Operators}

First, we collect some preliminary results on properties of the frame-related operators associated with localized sequences.

\begin{lem.}\label{lem:invertible}
Let $\varphi$ be a Riesz sequence in $\H$ that satisfies $\varphi\sim_\A \varphi$. Then $G_\varphi:\ell^{p}(X)\to\ell^p(X)$ is bounded and invertible for every $1\leq p\leq \infty$. 
\end{lem.}
\begin{proof}
By the localization assumption we have that $G_\varphi\in \A$. The Riesz sequence property implies that $G_\varphi$ is invertible on $\mathcal{B}(\ell^2(X))$, which implies by inverse closedness of $\A$ in $\mathcal{B}(\ell^2(X))$ that  $(G_\varphi)^{-1}\in \A$. Consequently, $G_\varphi$ is invertible on all $\ell^p(X)$-spaces ($1\leq p\leq \infty$). 
\end{proof}

\begin{cor.}\label{cor:D-phi-injective}
Let $\varphi$ be a Riesz basis for $\H$ that satisfies $\varphi\sim_\A \varphi$. Then $D_\varphi:\ell^p(X)\to \H^p(\varphi)$ is injective for every $1\leq p\leq \infty$.
\end{cor.}
\begin{proof}
We may write $G_\varphi:\ell^{p}(X)\to\ell^p(X)$ as $G_\varphi=C_\varphi D_\varphi$ with $D_\varphi:\ell^p(X)\to \H^p(\varphi)$ being bounded (Proposition~\ref{Donto}). The result then follows immediately from Lemma~\ref{lem:invertible}.
\end{proof}

The next two lemmata generalize \cite[Proposition~8]{groe04}, where the case of the Jaffard algebra $\A=\mathcal{J}_s$ is treated. They are also similar to \cite[Lemma~6 and  7]{xxlgro14}, except for the frame assumption on $\psi$. In fact, our conclusion is also weaker due to the missing frame assumption, i.e., we cannot conclude here that $C_\psi$ has closed range, as explicitly demonstrated in Example~\ref{counterexample} (iv). 

\begin{lem.}\label{lem:bessel}
Let $\varphi = (\varphi_k)_{k\in X}$ be a frame in $\Hil$ and $\psi = (\psi_k)_{k\in X}$ be a family in $\H$ that satisfies $\psi\sim_\A \varphi$. Then $C_\psi:\H^p(\varphi)\to \ell^p(X)$ is well-defined and bounded for each $1\leq p\leq \infty$. In particular, $\psi$ is a Bessel sequence in $\Hil$.
\end{lem.}

\begin{proof}
First, let $1\leq p <\infty$. Let $f\in \Hil^{00}(\varphi)$ be arbitrary. Then $C_{\widetilde{\varphi}}f \in \ell^p(X)$ and 
$$\Vert C_{\psi} f\Vert_{\ell^p(X)} = \Vert G_{\psi,\varphi} C_{\widetilde{\varphi}} f\Vert_{\ell^p(X)} \leq \Vert G_{\psi,\varphi} \Vert_{\B(\ell^p(X))} \Vert f\Vert_{\Hil^p(\varphi)},$$  
where we used frame reconstruction in $\Hil^{00}(\varphi)\subseteq \Hil$ with respect to $\varphi$ and $\widetilde{\varphi}$, as well as property (A1) from Definition~\ref{def:spectral}. Thus, by density, $C_{\psi}$ extends to a bounded operator $\H^p(\varphi)\to \ell^p(X)$ which we again denote by $C_{\psi}$.

Now let $p=\infty$. Let $f = \big[(  f_n )_{n\in\N}\big]_{\sim_{\varphi}}\in \Hil^{\infty}(\varphi)$ be arbitrary. Then, by definition, $\lim_{n\to\infty} \langle f_n, \widetilde{\varphi}_k \rangle =: \langle f,\widetilde{\varphi}_k \rangle$ converges for every $k\in X$ and $\sup_{n\in \mathbb{N}} \Vert C_{\widetilde{\varphi}} f_n\Vert_{\ell^{\infty}(X)}$ $\leq C$ for some $C>0$, as well as $C_{\widetilde{\varphi}}f := (\langle f,\widetilde{\varphi}_k\rangle)_{k\in X} \in \ell^{\infty}(X)$. This implies that 
\begin{equation}\label{eq:proof-lem16}\sum_{l\in X}\vert \langle f_n, \widetilde{\varphi}_l \rangle \vert \vert \langle\varphi_l, \psi_k \rangle\vert \leq C \sup_{k\in X}\sum_{l\in X} \vert \langle\varphi_l, \psi_k \rangle\vert = C \Vert G_{\psi,\varphi} \Vert_{\B(\ell^{\infty}(X))} <\infty 
\end{equation}
for each $k\in X$ and each $n\in \mathbb{N}$. Hence we may apply the dominated convergence theorem and obtain via frame reconstruction 
$$\lim_{n\to\infty} \langle f_n , \psi_k\rangle = \lim_{n\to\infty} \sum_{l\in X} \langle f_n, \widetilde{\varphi}_l \rangle \langle\varphi_l, \psi_k \rangle = \sum_{l\in X} \langle f, \widetilde{\varphi}_l \rangle \langle\varphi_l, \psi_k \rangle = [G_{\psi,\varphi} C_{\widetilde{\varphi}}f]_k$$
for each $k\in X$. Consequently, 
\begin{flalign}
\Vert C_{\psi} f\Vert_{\ell^{\infty}(X)} &= \left\Vert \left( \lim_{n\to\infty}\langle f_n , \psi_k\rangle\right)_{k\in X} \right\Vert_{\ell^{\infty}(X)} \notag \\
&= \Vert G_{\psi,\varphi} C_{\widetilde{\varphi}}f \Vert_{\ell^{\infty}(X)} \notag \\
&\leq \Vert G_{\psi,\varphi} \Vert_{\B(\ell^{\infty}(X))} \Vert C_{\widetilde{\varphi}}f \Vert_{\ell^{\infty}(X)} \notag \\
&= \Vert G_{\psi,\varphi} \Vert_{\B(\ell^{\infty}(X))} \Vert f\Vert_{\Hil^{\infty}(\varphi)}, \notag
\end{flalign}
as was to be shown.

Finally, boundedness of $C_{\psi}:\Hil^2(\varphi) \to \ell^2(X)$ is equivalent to $\psi$ being a Bessel sequence, since $\Hil^2(\varphi) = \Hil$.
\end{proof}

\begin{lem.}\label{Dbounded}
Let $\varphi = (\varphi_k)_{k\in X}$ be a frame in $\Hil$ and $\psi = (\psi_k)_{k\in X}$ be a family in $\H$ that satisfies $\psi\sim_\A \widetilde{\varphi}$. 
Then $D_\psi:\ell^p(X) \to \H^p(\varphi)$, defined analogously as in (\ref{syntheisseries}) and (\ref{synth}), is well-defined and bounded for each $1\leq p\leq \infty$ and the associated series converges unconditionally (resp., unconditionally in the $\sigma(\Hil,\Hil^{00}(\widetilde\varphi))$-topology in the case $p=\infty$).
\end{lem.}

\begin{proof}
Assume that $1\leq p <\infty$. Then, by \cite[Proposition~3.22]{koeba24}, $V^p(\varphi) = \lbrace f\in \Hil : C_{\widetilde{\varphi}}f \in \ell^p(X)\rbrace$ is a norm-dense subspace of $\Hil^p(\varphi)$. For $c=(c_l)_{l\in X} \in \ell^{00}(X)$ we have that $D_{\psi}c\in \Hil$ as it is a finite linear combination of elements from $\psi$. Since 
$$\Vert D_{\psi} c\Vert_{\Hil^p(\varphi)} = \Vert C_{\widetilde{\varphi}} D_{\psi} c\Vert_{\ell^p(X)} = \Vert G_{\widetilde{\varphi},\psi } c\Vert_{\ell^p(X)} \leq \Vert G_{\widetilde{\varphi},\psi }\Vert_{\B(\ell^p(X))} \Vert c\Vert_{\ell^p(X)},$$
we see that $D_{\psi}$ maps $\ell^{00}(X)$ boundedly into $V^p(\varphi)\subseteq \Hil^p(\varphi)$. Hence, the claim follows via density.

Now consider the case $p=\infty$. Pick some arbitrary $c=(c_l)_{l\in X} \in \ell^{\infty}(X)$ and assume that $F_1 \subseteq F_2\subseteq \dots$ is a nested sequence of finite subsets of $X$ such that $\bigcup_{n=1}^{\infty} F_n = X$. Then $f_n := \sum_{l\in F_n} c_l \psi_l$ is contained in $\Hil$ for each $n\in \mathbb{N}$ and 
\begin{equation}\label{compconv}
\lim_{n\rightarrow \infty} \langle f_n, \widetilde{\varphi}_k \rangle = \lim_{n\rightarrow \infty} \sum_{l\in F_n} \langle\psi_l, \widetilde{\varphi}_k \rangle c_l = \sum_{l\in X} [G_{\widetilde{\varphi},\psi }]_{k,l} c_l = [G_{\widetilde{\varphi},\psi } c]_k \in \mathbb{C}
\end{equation}
for each $k\in X$, since $G_{\widetilde{\varphi},\psi } \in \A \subseteq \B(\ell^{\infty}(X))$ by assumption. Furthermore, 
\begin{flalign}
\vert \langle f_n, \widetilde{\varphi}_k \rangle \vert &\leq \sum_{l\in F_n} \big\vert[G_{\widetilde{\varphi},\psi }]_{k,l} c_l \big\vert \notag \\
&\leq \sup_{k\in X}\sum_{l\in X} \big\vert[G_{\widetilde{\varphi},\psi }]_{k,l} \big\vert\,  \vert c_l \vert \notag  \\
&\leq \Big\Vert \big\vert G_{\widetilde{\varphi},\psi }\big\vert \Big\Vert_{\B(\ell^{\infty}(X))} \Vert c \Vert_{\ell^{\infty}(X)} \notag 
\end{flalign}
for all $k\in X$ and all $n\in \mathbb{N}$, where we used that $\big\vert G_{\widetilde{\varphi},\psi }\big\vert \in \A \subseteq \B(\ell^{\infty}(X))$ due to solidity of $\A$. This means that $D_{\psi}c = \big[( f_n )_{n\in\N}\big]_{\sim_{\varphi}} \in \Hil^{\infty}(\varphi)$. 
  The boundedness of $D_{\psi}$ follows immediately via (\ref{compconv}), as we have 
$$\Vert D_{\psi}c \Vert_{\Hil^{\infty}(\varphi)} = \sup_{k\in X} \big\vert [G_{\widetilde{\varphi},\psi } c]_k \big\vert \leq \Vert G_{\widetilde{\varphi},\psi }\Vert_{\B(\ell^{\infty}(X))} \Vert c \Vert_{\ell^{\infty}(X)}.$$
To show the unconditional convergence of $D_{\psi}c$ with respect to the $\sigma(\Hil, \Hil^{00}(\widetilde{\varphi}))$-topology, it suffices to show that $D_{\psi}c$ is independent of the choice of $\lbrace F_n \rbrace_{n=1}^{\infty}$ \cite[Section~3.2.4]{koeba24}. 
So, let 
$J_1 \subseteq J_2\subseteq \dots$ be another nested sequence of finite subsets of $X$ with $\bigcup_{n=1}^{\infty} J_n = X$. We will show that $\lbrace \sum_{l\in J_n}c_l \psi_l \rbrace_{n=1}^{\infty} \sim_{\varphi} \lbrace \sum_{l\in F_n}c_l \psi_l \rbrace_{n=1}^{\infty}$. Setting $H_n = F_n \cup J_n$ ($n\in \mathbb{N}$) gives a nested sequence of finite subsets whose union is $X$ as well. Then, 
for each $k\in X$,
\begin{flalign}
\lim_{n\rightarrow \infty} \left\vert \left\langle \sum_{l\in H_n} c_l \psi_l - \sum_{l\in F_n} c_l \psi_l , \widetilde{\varphi}_k \right\rangle \right\vert &= \lim_{n\rightarrow \infty} \left\vert \left\langle \sum_{l\in H_n\setminus F_n} c_l \psi_l , \widetilde{\varphi}_k \right\rangle \right\vert \notag \\
&\leq \Vert c \Vert_{\ell^{\infty}(X)}\lim_{n\rightarrow \infty} \sum_{l\in X\setminus F_n} \vert \langle \psi_l , \widetilde{\varphi}_k \rangle \vert  = 0,\notag
\end{flalign}
since 
$$\sup_{k\in X} \sum_{l\in X} \vert \langle \psi_l , \widetilde{\varphi}_k \rangle \vert = \Vert G_{\widetilde{\varphi},\psi} \Vert_{\B(\ell^{\infty}(X))} < \infty.$$
Clearly, it also holds 
$$\sup_{n\in \mathbb{N}}\left\Vert C_{\widetilde{\varphi}} \left( \sum_{l\in H_n} c_l \psi_l - \sum_{l\in F_n} c_l \psi_l \right) \right\Vert < \infty.$$
Thus we have shown that $\lbrace \sum_{l\in H_n}c_l \psi_l \rbrace_{n=1}^{\infty} \sim_{\varphi} \lbrace \sum_{l\in F_n}c_l \psi_l \rbrace_{n=1}^{\infty}$. The same argument holds if we replace $F_n$ by $J_n$. Thus, $\lbrace \sum_{l\in J_n}c_l \psi_l \rbrace_{n=1}^{\infty} \sim_{\varphi} \lbrace \sum_{l\in F_n}c_l \psi_l \rbrace_{n=1}^{\infty}$, since $\sim_{\varphi}$ is an equivalence relation.
\end{proof}

When we consider the operators $C_\psi:\H^p(\varphi)\to\ell^p(X)$ and $D_\psi:\ell^p(X)\to \H^p(\varphi)$, $1\leq p\leq \infty$, we write $\text{Ran}_p(C_\psi),\text{Ran}_p(D_\psi)$ and $\text{Ker}_p(C_\psi),\text{Ker}_p(D_\psi)$ for their respective ranges and nullspaces.

The next two lemmata are mild generalizations of parts of \cite[Lemma~7]{xxlgro14}.

\begin{lem.}\label{Dadjoint}
Let $\varphi = (\varphi_k)_{k\in X}$ be a frame in $\Hil$ satisfying $\varphi\sim_{\A}\varphi$ and $\psi = (\psi_k)_{k\in X}$ be a family in $\H$ that satisfies $\psi\sim_\A \varphi$. Let $1\leq p <\infty$ be arbitrary and $\frac{1}{p}+\frac{1}{q} =1$. Then the adjoint $D_{\psi}':\Hil^q(\varphi)\to \ell^q(X)$ of $D_\psi:\ell^p(X) \to \H^p(\varphi)$ is given by $C_{\psi}$.    
\end{lem.}

\begin{proof}
We first observe that for each $1\leq p <\infty$, we have $\psi_l\in \Hil^p(\varphi)$ for each $l\in X$ due to
$$\Vert \psi_l \Vert_{\Hil^p(\varphi)} = \Vert C_{\widetilde{\varphi}}\psi_l\Vert_{\ell^p(X)} \leq \Vert C_{\widetilde{\varphi}}\psi_l\Vert_{\ell^1(X)} \leq \sup_{l\in X}\sum_{k\in X}\vert \langle \psi_l, \widetilde{\varphi}_k\rangle \vert = \Vert G_{\widetilde{\varphi},\psi } \Vert_{\B(\ell^1(X))}.$$
We consider the case $(p,q)=(1,\infty)$. Let $f=\big[( f_n )_{n\in\N}\big]_{\sim_{\varphi}}\in \H^{\infty}(\varphi)$ and $c=(c_k)_{k\in X}\in \ell^1(X)$ be arbitrary. Then
$$\langle c,D_{\psi}'f \rangle_{\ell^1 \times \ell^{\infty}} = \langle D_{\psi}c,f\rangle_{\Hil^1 \times \Hil^{\infty}} = \sum_{k\in X}c_k \langle \psi_k,f\rangle_{\Hil^1 \times \Hil^{\infty}}.$$
According to Theorem~\ref{duality}, we have that
\begin{flalign}
\langle \psi_k,f\rangle_{\Hil^1 \times \Hil^{\infty}} &= \sum_{l\in X} \langle \psi_k, \widetilde{\varphi}_l\rangle \lim_{n\to \infty}\langle \varphi_l, f_n\rangle \notag \\ 
&= \lim_{n\to\infty}\sum_{l\in X} \langle \psi_k, \widetilde{\varphi}_l\rangle\langle \varphi_l, f_n\rangle \notag \\
&= \lim_{n\to\infty} \langle \psi_k, f_n\rangle = \overline{[C_{\psi} f]_k}, \notag
\end{flalign}
where we applied the dominated convergence theorem (justified analogously as in \eqref{eq:proof-lem16} 
upon noticing that $\Vert C_{\varphi} \cdot \Vert_{\ell^{\infty}(X)} \asymp \Vert C_{\widetilde{\varphi}} \cdot \Vert_{\ell^{\infty}(X)}$ (c.f. Proposition~\ref{Hinftynormequ})) in the second step and frame reconstruction in the third step. On the other hand, 
$$\langle c,C_{\psi}f \rangle_{\ell^1 \times \ell^{\infty}} = \sum_{k\in X}c_k \overline{[C_{\psi} f]_k},$$
which yields the claim. For the other values of $p$ and $q$ we proceed analogously by showing the identity for $f\in \H^{00}(\varphi)$ and then argue via density. 
\end{proof}

\begin{lem.}\label{Cadjoint}
Let $\varphi = (\varphi_k)_{k\in X}$ be a frame in $\Hil$ satisfying $\varphi\sim_{\A}\varphi$ and $\psi = (\psi_k)_{k\in X}$ be a family in $\H$ that satisfies $\psi\sim_\A \varphi$. Let $1\leq p <\infty$ be arbitrary and $\frac{1}{p}+\frac{1}{q} =1$. Then the Banach space adjoint $C_{\psi}':\ell^q(X) \to \H^q(\varphi)$ of $C_\psi:\Hil^p(\varphi)\to \ell^p(X)$ is given by $D_{\psi}$.
\end{lem.}

\begin{proof}
We first consider the case $(p,q)=(1,\infty)$. Let $f\in \H^{00}(\varphi) \subseteq \H^1(\varphi)$ and $c=(c_l)_{l\in X}\in \ell^{\infty}$ be arbitrary. According to Lemma~\ref{Dbounded}, $D_{\psi}c =\big[( f_n)_{n\in\N}\big]_{\sim_{\varphi}} \in \H^{\infty}(\varphi)$, where $f_n = \sum_{l\in F_n}c_l \psi_l$ and $F_1 \subseteq F_2\subseteq \dots$ is a nested sequence of finite subsets of $X$ such that $\bigcup_{n=1}^{\infty} F_n = X$. Then we see that 
\begin{flalign}
\langle f, C_{\psi}'c\rangle_{\Hil^1\times \Hil^{\infty}} &= \langle C_{\psi}f,c \rangle_{\ell^1\times\ell^{\infty}} \notag \\    
&= \sum_{l\in X} \langle f,\psi_l\rangle \overline{c_l}\notag \\
&= \lim_{n\to\infty} \left\langle f, \sum_{l\in F_n} c_l \psi_l \right\rangle \notag \\
&= \lim_{n\to\infty} \left\langle \sum_{k\in X} \langle f,\widetilde{\varphi}_k\rangle \varphi_k , \sum_{l\in F_n} c_l \psi_l \right\rangle \notag \\
&= \sum_{k\in X} \langle f,\widetilde{\varphi}_k\rangle \lim_{n\to\infty} \left\langle \varphi_k , \sum_{l\in F_n} c_l \psi_l \right\rangle \notag \\
&= \langle f, D_{\psi}c \rangle_{\Hil^1\times \Hil^{\infty}}, \notag 
\end{flalign}
where we used frame reconstruction in the fourth step and the dominated convergence theorem (justified analogously as in \eqref{eq:proof-lem16} upon noticing that $\Vert C_{\varphi} \cdot \Vert_{\ell^{\infty}(X)} \asymp \Vert C_{\widetilde{\varphi}} \cdot \Vert_{\ell^{\infty}(X)}$ (c.f. Proposition~\ref{Hinftynormequ})) in the fifth step. Via density we conclude that $\langle f, C_{\psi}'c\rangle_{\Hil^1\times \Hil^{\infty}} = \langle f, D_{\psi}c \rangle_{\Hil^1\times \Hil^{\infty}}$ for all $f\in \Hil^1(\varphi)$ and the claim is proven. For other values of $p$ and $q$ the proof is similar.
\end{proof}

As a consequence we immediately obtain the following.

\begin{cor.}\label{SGadoints}
Let $\varphi = (\varphi_k)_{k\in X}$ be a frame in $\Hil$ satisfying $\varphi\sim_{\A}\varphi$ and $\psi = (\psi_k)_{k\in X}$ be a family in $\H$ that satisfies $\psi\sim_\A \varphi$. For $1\leq p <\infty$ and $\frac{1}{p}+\frac{1}{q} = 1$ the following hold. 
\begin{itemize}
\item[(i)] $S_{\psi} \in \B(\H^p(\varphi))$ and $S_{\psi}' = S_{\psi}\in \B(\H^q(\varphi))$.
\item[(ii)] $G_{\psi} \in \B(\ell^p(X))$ and $G_{\psi}' = G_{\psi} \in \B(\ell^q(X))$.
\end{itemize}
\end{cor.}

The following consequence of Remark~\ref{holfunccalc} on the holomorphic functional calculus is fundamental for the proof of our main result.

\begin{lem.}\label{alpha}
Let $\varphi$ be an intrinsically $\A$-localized Riesz basis in $\Hil$. Then $G_{\varphi}^{\alpha} \in \A$ for every $\alpha \in \mathbb{R}$.
\end{lem.}

\begin{proof}
Since $\A$ is inverse-closed in $\B(\ell^2(X))$, the spectral identity $ \sigma_{\B(\ell^2(X))}(A)=\sigma_{\mathcal{A}}(A)$ holds true for all $A\in \A$, and the holomorphic functional calculi associated with $\A$ and $\B(\ell^2(X))$ coincide (see Remark~\ref{holfunccalc}). In particular, if $A\leq B$ denote the Riesz basis constants of $\varphi$, then $\sigma_{\A}(G_{\varphi}) = \sigma_{\B(\ell^2(X))}(G_{\varphi}) \subseteq [A,B]$. Thus, there exists a closed path $\Gamma$ contained in the right half plane $\text{Re}(z)>0$ with $\sigma_{\B(\ell^2(X))}(G_{\varphi}) = \sigma_{\A}(G_{\varphi}) \subset \text{int}(\Gamma)$. Since $f(z) = z^{\alpha}$ is holomorphic on $\text{Re}(z)>0$, the holomorphic functional calculus yields $G_{\varphi}^{\alpha}\in \A$.  
\end{proof}

\begin{cor.}\label{easycorollary}
Let $\varphi$ be an intrinsically $\A$-localized Riesz basis in $\Hil$. Then both $G_{S_{\varphi}^{-1/4}\varphi}$ and $G_{S_{\varphi}^{-3/4}\varphi}$ are contained in $\A$.
\end{cor.}

\begin{proof}
Since $\varphi$ is a Riesz basis, $(S_{\varphi}^{-1/2}\varphi_k)_{k\in X}$ is an orthonormal basis. By a simple calculation we see that $\big(G_{S_{\varphi}^{-1/4}\varphi}\big)^2 = G_{\varphi}$. Thus, by Lemma~\ref{alpha}, $G_{S_{\varphi}^{-1/4}\varphi} = G_{\varphi}^{1/2}\in \A$. Similarly, we have $(G_{S_\varphi^{-3/4}\varphi})^2=G_{\widetilde\varphi} \in \A$ due to Theorem~\ref{dualloc}, which, by Lemma~\ref{alpha}, implies $G_{S_\varphi^{-3/4}\varphi}\in\mathcal{A}$.
\end{proof}

\begin{lem.}\label{powers-of-S}
Let $\varphi$ be an intrinsically $\A$-localized Riesz basis. Then both $S_{\varphi}^{1/2}$ and $S_{\varphi}^{-1/2}$ are well-defined and bounded on $\Hil^p(\varphi)$ for all $1\leq p \leq \infty$. Furthermore, both $\big(S_{\varphi}^{1/2}\big)' = S_{\varphi}^{1/2}$ and $\big(S_{\varphi}^{-1/2}\big)' = S_{\varphi}^{-1/2}$ are self-adjoint.
\end{lem.}

\begin{proof}
Let $\alpha \in \lbrace - \frac{1}{2}, \frac{1}{2} \rbrace$. Then, by Lemma \ref{alpha}, $G_{\varphi}^{\alpha} \in \A$, so that $D_{\varphi} G_{\varphi}^{\alpha} C_{\widetilde{\varphi}} \in \B(\Hil^p(\varphi))$ for all $1\leq p \leq \infty$. Next, by \cite[Proposition~2.30]{hol22}, both $D_{\varphi}:\ell^p(X) \to \Hil^p(\varphi)$ and $C_{\varphi}:\Hil^p(\varphi) \to \ell^p(X)$ are invertible for all $1\leq p \leq \infty$ (alternatively, this also follows from Corollary \ref{cor:D-phi-injective}). In combination with the fact that $S_{\varphi} \in \B(\Hil^p(\varphi))$ being invertible (Theorem \ref{Sinvertible}) and $C_{\widetilde{\varphi}} = C_{\varphi} S_{\varphi}^{-1}$, we see
$$\big(D_{\varphi} G_{\varphi}^{-1/2} C_{\widetilde{\varphi}}\big)^2 = D_{\varphi} G_{\varphi}^{-1} C_{\varphi} S_{\varphi}^{-1} = D_{\varphi} D_{\varphi}^{-1} C_{\varphi}^{-1} C_{\varphi} S_{\varphi}^{-1} = S_{\varphi}^{-1}$$
and similarly $\big(D_{\varphi} G_{\varphi}^{1/2} C_{\widetilde{\varphi}}\big)^2 = S_{\varphi}$, so that $S_{\varphi}^{\alpha} = D_{\varphi} G_{\varphi}^{\alpha} C_{\widetilde{\varphi}}$. 

Now let $1\leq p <\infty$ and $\frac{1}{p}+\frac{1}{q} =1$. Using $C_{\varphi} = G_{\varphi} C_{\widetilde{\varphi}}$ and the self-adjointness of $G_{\varphi}^{-1/2} = G_{S_{\varphi}^{-1/4} \varphi}$ (see Corollary \ref{easycorollary}), we see that for all $f\in \Hil^p(\varphi), g\in \Hil^{q}(\varphi)$
\begin{flalign}
    \langle S_{\varphi}^{-1/2} f,g \rangle_{\Hil^p \times \Hil^{q}} &= \langle C_{\widetilde{\varphi}} S_{\varphi}^{-1/2} f, C_{\varphi} g \rangle_{\ell^p \times \ell^{q}} \notag \\
    &= \langle C_{\widetilde{\varphi}} D_{\varphi} G_{\varphi}^{-1/2} C_{\widetilde{\varphi}} f, C_{\varphi} g \rangle_{\ell^p \times \ell^{q}} \notag \\
    &= \langle G_{\varphi}^{-1/2} C_{\widetilde{\varphi}} f, C_{\varphi} g \rangle_{\ell^p \times \ell^{q}} \notag \\
    &= \langle  C_{\widetilde{\varphi}} f, G_{\varphi}^{-1/2} G_{\varphi} C_{\widetilde{\varphi}} g \rangle_{\ell^p \times \ell^{q}} \notag \\
    &= \langle  C_{\widetilde{\varphi}} f, C_{\varphi} D_{\varphi} G_{\varphi}^{-1/2}C_{\widetilde{\varphi}} g \rangle_{\ell^p \times \ell^{q}} \notag \\
    &= \langle  C_{\widetilde{\varphi}} f, C_{\varphi} S_{\varphi}^{-1/2} g \rangle_{\ell^p \times \ell^{q}}  \notag \\
    &= \langle  f, S_{\varphi}^{-1/2} g \rangle_{\Hil^p \times \Hil^{q}} .  \notag 
\end{flalign}
The corresponding duality relation for $S_{\varphi}^{1/2}$ can be shown similarly.
\end{proof}

\subsection{Localization of the R-dual sequence}

One main ingredient in Theorem~\ref{thm:main} is the family $\omega$, which, in the language of \cite{r-duals}, is a so-called \emph{R-dual of type IV} of $\psi$. 

Let $(\psi_k)_{k\in X}$ be a Bessel sequence in $\Hil$ and let $\varphi = (\varphi_k)_{k\in X}$ be a Riesz basis in $\Hil$. Then $(S_{\varphi}^{-1/2}\varphi_k)_{k\in X}$ is an orthonormal basis in $\Hil$ \cite{ole1n}. Consequently, 
$$\omega_k = \sum_{l\in X}\langle \psi_l , \varphi_k \rangle S_{\varphi}^{- {1}/{2}}\varphi_l$$
converges in $\Hil$ for each $k\in X$, defining a family $\omega = (\omega_k)_{k\in X}$ of vectors in $\Hil$. One important observation for proving Theorem~\ref{thm:main} is, that 
$\omega \sim_{\A}\varphi$ whenever $\psi \sim_{\A}\varphi$, see below. 

\begin{lem.}\label{omegalocalization}
Let $\A$ be a spectral algebra and $\varphi = (\varphi_k)_{k\in X}$ be a Riesz basis in $\Hil$ that satisfies $\varphi\sim_{\A}\varphi$. Let $\psi = (\psi_k)_{k\in X}$ be a sequence in $\Hil$ satisfying $\psi \sim_{\A} \varphi$ and let $\omega = (\omega_k)_{k\in X}$ be given by 
$$\omega_k = \sum_{l\in X}\langle \psi_l , \varphi_k \rangle S_{\varphi}^{- {1}/{2}}\varphi_l \qquad (k\in X).$$
Then $\omega \sim_{\A}\varphi$, $\omega \sim_{\A}\widetilde{\varphi}$ and $\omega \sim_{\A}\omega$. Consequently, $\omega$ is a Bessel sequence in $\Hil$.
\end{lem.}

\begin{proof}
First, note that $\psi$ is a Bessel sequence, since $G_\psi \in \A \subset \B(\ell^2(X))$ according to Lemma~\ref{easy}, whence $\omega$ is a well-defined sequence in $\Hil$. Next, recall from Corollary~\ref{easycorollary} that $G_{S_{\varphi}^{-1/4}\varphi} \in \A$. Since 
$$\langle \varphi_n, \omega_k\rangle = \sum_{l\in X}\langle \varphi_k,\psi_l\rangle \left\langle S_{\varphi}^{-1/4}\varphi_n, S_{\varphi}^{-1/4} \varphi_l\right\rangle = \sum_{l\in X}[G_{\psi,\varphi}^*]_{k,l} \left[G_{S_{\varphi}^{-1/4}\varphi}\right]_{l,n}$$
for all $k,n\in X$, we see that $G_{\omega,\varphi} = G_{\psi,\varphi}^* G_{S_{\varphi}^{-1/4}\varphi} \in \A$.
 
Similarly, since $G_{S_\varphi^{-3/4}\varphi}\in\mathcal{A}$ and
$$\langle \widetilde{\varphi}_n, \omega_k\rangle = \sum_{l\in X}\langle \varphi_k, \psi_l\rangle \left\langle S_{\varphi}^{-3/4}\varphi_n, S_{\varphi}^{-3/4} \varphi_l\right\rangle = \sum_{l\in X}[G_{\psi,\varphi}^*]_{k,l} \left[G_{S_{\varphi}^{-3/4}\varphi}\right]_{l,n}$$
for all $k,n\in X$, we obtain $G_{\omega,\widetilde{\varphi}} = G_{\psi,\varphi}^* G_{S_{\varphi}^{-3/4}\varphi} \in \A$. 

Finally, $G_{\omega} = G_{\omega,\widetilde{\varphi}}G_{{\varphi},\omega} \in \A$, hence $G_{\omega}\in \B(\ell^2(X))$, which implies that $\omega$ is a Bessel sequence in $\Hil$.
\end{proof}

As a consequence, all results appearing in the previous subsection apply to $\omega$, which will be crucial for proving Theorem~\ref{thm:main}. 

Finally, we define the operators 
$$\overline{C_\psi}f=\big(\langle \psi_k,f\rangle\big)_{k\in X} \qquad \text{and} \qquad \overline{D_\psi} c = \sum_{l\in X} \overline{c_l} \psi_l = D_{\psi} \overline{c}$$ to rewrite $D_\omega$ and $C_\omega$ in terms of $\psi$ and $\varphi$.

\begin{lem.}\label{rewrite:D_wc}
Let $\psi,\varphi$ and $\omega$ be as in Lemma~\ref{omegalocalization}. Then, for all $1\leq p \leq \infty$, the operators $D_\omega:\ell^p(X)\to \H^p(\varphi)$ and $C_\omega: \H^p(\varphi) \to \ell^p(X)$ equal 
$$D_\omega =S_{\varphi}^{-1/2} D_{\varphi} \overline{C_{\psi}} \, \overline{D_{\varphi}} \qquad \text{and} \qquad C_\omega = \overline{C_{\varphi}} \,\overline{D_{\psi}} C_\varphi S_{\varphi}^{-1/2}.$$
\end{lem.}

\begin{proof} By our previous results, all occurring operators are well-defined and bounded on their respective domains. Furthermore, by duality (using Lemmata \ref{Dadjoint}, \ref{Cadjoint} and \ref{powers-of-S}) it suffices to show the claimed operator identities for $1\leq p \leq 2$, and by density, it suffices to verify them on $\ell^{00}(X)$, respectively on $\Hil^{00}(\varphi)$. To this end, we have for all $c\in \ell^{00}(X)$
\begin{flalign}
    D_\omega c &= \sum_{k\in X} c_k \sum_{l\in X}\langle \psi_l , \varphi_k \rangle S_{\varphi}^{- {1}/{2}}\varphi_l \notag \\
    &= \sum_{l\in X} \sum_{k\in X}\langle \psi_l , \overline{c_k} \varphi_k \rangle S_{\varphi}^{- {1}/{2}}\varphi_l \notag \\
    &= S_{\varphi}^{-1/2} \sum_{l\in X} \langle \psi_l , \overline{D_{\varphi}} c \rangle \varphi_l = S_{\varphi}^{-1/2} D_{\varphi} \overline{C_{\psi}} \, \overline{D_{\varphi}} c .\notag
\end{flalign}
Likewise, for all $f\in \Hil^{00}(\varphi)$, we have
\begin{flalign}
C_\omega f &= \bigg( \Big\langle f, \sum_{l\in X}\langle \psi_l , \varphi_k \rangle S_{\varphi}^{- {1}/{2}}\varphi_l \Big\rangle \bigg)_{k\in X} \notag \\
&= \bigg( \sum_{l\in X} \langle \varphi_k , \psi_l \rangle \langle  f, S_{\varphi}^{-1/2} \varphi_l \rangle \bigg)_{k\in X} \notag \\
&= \bigg( \Big\langle \varphi_k , \sum_{l\in X} \overline{\big[ C_\varphi S_{\varphi}^{-1/2} f\big]_l} \psi_l \Big\rangle \bigg)_{k\in X} = \overline{C_{\varphi}} \,\overline{D_{\psi}} C_\varphi S_{\varphi}^{-1/2} f , \notag
\end{flalign}
as was to be shown.
%
%
\end{proof}

\section{Proof of Theorem~\ref{thm:main}}\label{Proof of main}

\begin{figure}\begin{center}
\begin{tikzpicture}[scale=0.5]
    \draw[thick] ({4*cos(360/6+30)},{4*sin(360/6+30)}) circle [radius=0.4];
    \node[thick] at ({4*cos(360/6+30)},{4*sin(360/6+30)}) {\textbf{\tiny{1}}};

       \draw[thick] ({4*cos(2*360/6+30)},{4*sin(2*360/6+30)}) circle [radius=0.4];
    \node[thick] at ({4*cos(2*360/6+30)},{4*sin(2*360/6+30)}) {\textbf{\tiny{2}}};

         \draw[thick] ({6.5*cos(2*360/6+30)},{6.5*sin(2*360/6+30)}) circle [radius=0.4];
    \node[thick] at ({6.5*cos(2*360/6+30)},{6.5*sin(2*360/6+30)}) {\textbf{\tiny{3}}};

       \draw[thick] ({4*cos(3*360/6+30)},{4*sin(3*360/6+30)}) circle [radius=0.4];
    \node[thick] at ({4*cos(3*360/6+30)},{4*sin(3*360/6+30)}) {\textbf{\tiny{5}}};

       \draw[thick] ({4*cos(4*360/6+30)},{4*sin(4*360/6+30)}) circle [radius=0.4];
    \node at ({4*cos(4*360/6+30)},{4*sin(4*360/6+30)}) {\textbf{\tiny{6}}};

       \draw[thick] ({4*cos(5*360/6+30)},{4*sin(5*360/6+30)}) circle [radius=0.4];
    \node[thick] at ({4*cos(5*360/6+30)},{4*sin(5*360/6+30)}) {\textbf{\tiny{7}}};

           \draw[thick] ({6.5*cos(3*360/6+30)},{6.5*sin(3*360/6+30)}) circle [radius=0.4];
    \node[thick] at ({6.5*cos(3*360/6+30)},{6.5*sin(3*360/6+30)}) {\textbf{\tiny{4}}};


       \draw[thick] ({4*cos(6*360/6+30)},{4*sin(6*360/6+30)}) circle [radius=0.4];
    \node[thick] at ({4*cos(6*360/6+30)},{4*sin(6*360/6+30)}) {\textbf{\tiny{10}}};

          \draw[thick] ({6.5*cos(5*360/6+30)},{6.5*sin(5*360/6+30)}) circle [radius=0.4];
    \node[thick] at ({6.5*cos(5*360/6+30)},{6.5*sin(5*360/6+30)}) {\textbf{\tiny{8}}};

        \draw[thick] ({6.5*cos(6*360/6+30)},{6.5*sin(6*360/6+30)}) circle [radius=0.4];
    \node[thick] at ({6.5*cos(6*360/6+30)},{6.5*sin(6*360/6+30)}) {\textbf{\tiny{9}}};






 \draw[->,thick]  ({4.45*cos(360+30)},{4.45*sin(360+30)})--({6.05*cos(360+30)},{6.05*sin(360+30)});

  \draw[<-,thick]  ({4.45*cos(300+30)},{4.45*sin(300+30)})--({6.05*cos(300+30)},{6.05*sin(300+30)});

\draw[<-,thick,domain=305:355] plot ({6.5*cos(\x+30)},{6.5*sin(\x+30)});

    \draw[<->,thick,domain=7:53] plot ({4*cos(\x+30)},{4*sin(\x+30)});

      \draw[->,thick,domain=67:113] plot ({4*cos(\x+30)},{4*sin(\x+30)});

        \draw[->,thick,domain=127:173] plot ({4*cos(\x+30)},{4*sin(\x+30)});

        \draw[->,thick,domain=187:233] plot ({4*cos(\x+30)},{4*sin(\x+30)});

        \draw[<->,thick]  ({4.45*cos(120+30)},{4.45*sin(120+30)})--({6.05*cos(120+30)},{6.05*sin(120+30)});

 \draw[<->,thick]  ({4.45*cos(180+30)},{4.45*sin(180+30)})--({6.05*cos(180+30)},{6.05*sin(180+30)});

          \draw[<->,thick,domain=247:293] plot ({4*cos(\x+30)},{4*sin(\x+30)});

            \draw[->,thick,domain=307:353] plot ({4*cos(\x+30)},{4*sin(\x+30)});

\end{tikzpicture}
\end{center}
    \caption{Graph showing which implications of the equivalent statements in Theorem~\ref{thm:main} are shown.}
    \label{fig:picto1}
\end{figure}
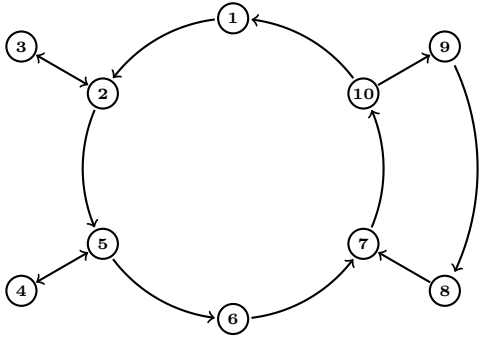

\begin{proof} We show the equivalences of Theorem~\ref{thm:main} via a chain of implications that is depicted in Figure~\ref{fig:picto1}.
Before we start, observe that by Lemma~\ref{omegalocalization} as well as Lemma~\ref{lem:bessel} and Lemma~\ref{Dbounded}, all operators occurring in the statements \textbf{\ref{list:2}}$-$\textbf{\ref{list:9}} are bounded on their respective domains.\medskip

\textbf{\ref{list:1}$\Rightarrow$\ref{list:2}:} The assumptions imply that $\psi$ is an intrinsically $\A$-localized frame, since $G_{\psi} = G_{\psi,\varphi} G_{ \widetilde{\varphi}} G_{\varphi,\psi}  \in \A$. Consequently, by Theorem~\ref{Sinvertible}, $S_{\psi}$ is bounded and invertible on $\Hil^1(\psi)$. By Remark~\ref{indep} it holds that $\Hil^1(\psi) = \Hil^1(\varphi)$ with equivalent norms. Hence (2) follows.  
\medskip

\textbf{\ref{list:2}$\Leftrightarrow$\ref{list:3}:} This follows from duality via Corollary~\ref{SGadoints} (i).\medskip


\medskip

\textbf{\ref{list:2}$\Rightarrow$\ref{list:5}:} Since $S_\psi=D_\psi C_\psi : \H^1(\varphi) \to \H^1(\varphi)$ is bijective, $D_\psi:\ell^1(X)\to \H^1(\varphi)$ has to be surjective. Furthermore, for all $f\in \Hil^1(\varphi)$ it holds 
$$\Vert f \Vert_{\Hil^1(\varphi)} = \Vert S_\psi^{-1} D_{\psi} C_{\psi} f \Vert_{\Hil^1(\varphi)} \leq \Vert S_\psi^{-1} D_{\psi} \Vert_{\B(\ell^1(X),\Hil^1(\varphi))} \Vert C_{\psi} f \Vert_{\ell^1(X)}.$$
This means that $C_{\psi}:\Hil^1(\varphi) \to \ell^1(X)$ is bounded from below, or, equivalently, injective and has closed range.
\medskip

\textbf{\ref{list:4}$\Leftrightarrow$\ref{list:5}:} By duality, Lemma~\ref{Dadjoint} and the closed range theorem, $C_{\psi}:\Hil^{\infty}(\varphi) \to \ell^{\infty}(X)$ is injective and has closed range if and only if $D_{\psi}:\ell^1(X)\to \Hil^1(\varphi)$ is surjective. By analogous reasoning via Lemma~\ref{Cadjoint}, $\text{Ran}_{1}(C_\psi)$ is closed if and only if $\text{Ran}_{\infty}(D_\psi)$ is closed.\medskip

\textbf{\ref{list:5}$\Rightarrow$\ref{list:6}:} First, we show that $D_{\omega}: \ell^{\infty}\to \Hil^{\infty}(\varphi)$ is injective and has closed range, or equivalently, that $D_\omega$ is bounded from below. 
By Lemma \ref{rewrite:D_wc}, $D_\omega c = S_{\varphi}^{-1/2} D_{\varphi} \overline{C_{\psi}} D_{\varphi} \overline{c}$ for all $c\in \ell^{\infty}(X)$. As observed in the proof of Lemma \ref{rewrite:D_wc}, both $S_{\varphi}^{-1/2}$ and  $D_{\varphi}$ are bounded and invertible, hence also bounded from below. Note that also $ \overline{C_{\psi}}$ is bounded from below, because $C_{\psi}$ is injective and has closed range due to \textbf{\ref{list:4}$\Leftrightarrow$\ref{list:5}}. Therefore $D_\omega$ is bounded from below, since $\Vert \overline{c} \Vert_{\ell^{\infty}} = \Vert c \Vert_{\ell^{\infty}}$.

Now, we show that $\text{Ran}_\infty(C_\omega)$ is closed. By Lemma~\ref{Dadjoint} and the closed range theorem, it suffices to prove that $D_\omega : \ell^1(X)\to \Hil^1(\varphi)$ has closed range. To this end, let $( D_{\omega} c^{(m)})_{m\in\N} \subset \Hil^{1}(\varphi)$ be a sequence in $\text{Ran}_1(D_{\omega})$ which converges in $\Hil^{1}(\varphi)$ to some $f \in \Hil^{1}(\varphi)$. By Lemma~\ref{rewrite:D_wc}, 
\begin{flalign}
0 &= \lim_{m\to\infty} \left\Vert D_{\omega} c^{(m)} -f \right\Vert_{\Hil^{1}(\varphi)} = \lim_{m\to\infty}\left\Vert C_{\widetilde{\varphi}} S_{\varphi}^{-1/2} D_{\varphi} \overline{C_{\psi}}D_{\varphi} \overline{c^{(m)}} - C_{\widetilde{\varphi}}f \right\Vert_{\ell^{1}(X)}. \notag
\end{flalign}
By continuity of $C_{\widetilde{\varphi}} S_{\varphi}^{1/2} D_{\varphi} \in \B(\ell^{1}(X))$, the latter implies 
$$0 = \lim_{m\to\infty}\left\Vert \overline{C_{\psi}}D_{\varphi} \overline{c^{(m)}} - C_{\widetilde{\varphi}}S_{\varphi}^{1/2}f \right\Vert_{\ell^{1}(X)},$$
which means that $\big( \overline{C_{\psi}}D_{\varphi} \overline{c^{(m)}} \big)_{m\in\N}$ is a sequence in $\text{Ran}_1\left(\overline{C_{\psi}}\right)$ which converges in $\ell^{1}(X)$ to $C_{\widetilde{\varphi}}S_{\varphi}^{1/2}f\in \ell^{1}(X)$. Since $\text{Ran}_1\left(\overline{C_{\psi}}\right)$ is closed due to $\text{Ran}_1\left(C_{\psi}\right)$ being closed, we must have $C_{\widetilde{\varphi}}S_{\varphi}^{1/2}f = \overline{C_{\psi}} h$ for some $h\in \Hil^{1}(\varphi)$. Since $D_{\varphi}:\ell^{1}(X)\to \Hil^{1}(\varphi)$ is surjective due to Proposition~\ref{Donto}, there exists some $\overline{c}\in \ell^{1}(X)$ with $h= D_{\varphi} \overline{c}$. Altogether, we obtain 
$$0 = \lim_{m\to\infty} \left\Vert \overline{C_{\psi}}D_{\varphi} \overline{c^{(m)}} - \overline{C_{\psi}} D_{\varphi} \overline{c} \right\Vert_{\ell^{1}(X)}.$$
Since $C_{\widetilde{\varphi}} S_{\varphi}^{-1/2} D_{\varphi} \in \B(\ell^{1}(X))$, the latter yields
$$ 0 = \lim_{m\to\infty}\left\Vert C_{\widetilde{\varphi}} S_{\varphi}^{-1/2}\ D_{\varphi} \overline{C_{\psi}}D_{\varphi} (\overline{c^{(m)}} - \overline{c}) \right\Vert_{\ell^{1}(X)}.$$
Lemma \ref{rewrite:D_wc} implies 
$$0 = \lim_{m\to\infty} \left\Vert D_{\omega} c^{(m)} -D_{\omega} c \right\Vert_{\Hil^{1}(\varphi)},$$
as was to be shown.
\medskip

\textbf{\ref{list:6}$\Leftrightarrow$\ref{list:7}:}  Lemma~\ref{omegalocalization} ensures that we may apply Lemma~\ref{Cadjoint} to deduce that the Banach space adjoint of $C_{\omega}:\Hil^1(\varphi)\to \ell^1(X)$ is $D_{\omega}:\ell^\infty(X)\to\Hil^\infty(\varphi)$. Therefore, $\text{Ker}_\infty(D_\omega)=\{0\}$ if and only if $\text{Ran}_1(C_\omega)$ is dense in $\ell^1(X)$.   The equivalence of the conditions on the closedness of $\text{Ran}_\infty(C_\omega)$ and $\text{Ran}_1(D_\omega)$ as well as $\text{Ran}_\infty(D_\omega)$ and $\text{Ran}_1(C_\omega)$ is due to the closed range theorem.

\medskip

\textbf{\ref{list:7}$\Rightarrow$\ref{list:10}:} By Lemma \ref{omegalocalization}, $\omega$ is a Bessel sequence, so that $D_\omega$ is well-defined and bounded on $\ell^2(X)$. Hence, it suffices to show that $D_\omega : \ell^2(X) \to \Hil$ is bounded from below. Since $\text{Ran}_1(C_\omega)=\ell^1(X)$, it follows that $\text{Ran}_1(C_\omega)$ is dense in $\ell^2(X)$ and so is $\text{Ran}_2(C_\omega)$. We therefore deduce that $\text{Ker}_1(D_\omega)\subset \text{Ker}_2(D_\omega)=\{0\}$, that is, $D_\omega:\ell^1(X)\to \Hil^1(\varphi)$ is injective. Since $\text{Ran}_1(D_\omega)$ is closed, it follows from the bounded inverse theorem that $D_\omega:\ell^1(X) \to \Hil^1(\varphi)$ is bounded from below. In particular, there exists \(A_1>0\) such that
$$A_1\|c\|_{\ell^1(X)} \leq \|D_\omega c\|_{\mathcal H^1(\varphi)} = \Vert C_{\widetilde{\varphi}} D_\omega c\Vert_{\ell^1(X)} = \Vert G_{\widetilde{\varphi}, \omega} c\Vert_{\ell^1(X)}, \qquad c\in\ell^1(X).$$
Thus, $G_{\widetilde{\varphi}, \omega}$ is bounded from below on $\ell^1(X)$ and, in the notation of Definition \ref{def:p-stability}, it holds \(\lambda_1(G_{\widetilde{\varphi}, \omega}) \geq A_1 > 0\). By Lemma \ref{omegalocalization}, $G_{\widetilde{\varphi}, \omega} \in \A$, and since $\A$ is uniformly $p$-stable by assumption, we conclude $\lambda_2(G_{\widetilde{\varphi}, \omega}) > 0$. Thus, there exists $A_2 > 0$ such that
$$A_2\|c\|_{\ell^2(X)} \leq \Vert G_{\widetilde{\varphi}, \omega} c\Vert_{\ell^2(X)} = \Vert C_{\widetilde{\varphi}} D_\omega c\Vert_{\ell^2(X)} = \|D_\omega c\|_{\mathcal H^2(\varphi)}, \qquad c\in\ell^2(X).$$ Since $\mathcal H^2(\varphi) = \Hil$ with equivalent norms, this means that $D_\omega : \ell^2(X) \to \Hil$ is bounded from below. \medskip

\textbf{\ref{list:10}$\Leftrightarrow$\ref{list:1}:} Note that $\big(S_\varphi^{1/2}\psi_k\big)_{k\in X}$ is a Bessel sequence because $(\psi_k)_{k\in X}$ is a Bessel sequence by Lemma~\ref{lem:bessel} and $S_\varphi^{1/2}\in \B(\Hil)$. Since $S_\varphi^{1/2}$ is invertible (see, e.g., \cite{ole1n}), we may rewrite $\omega_k$ as
$$\omega_k=\sum_{l\in X}\langle S_\varphi^{1/2}\psi_l,S_\varphi^{-1/2}\varphi_k\rangle S_\varphi^{-1/2}\varphi_l,$$
which means that $(\omega_k)_{k\in X}$ is an \emph{R-dual} sequence of $\big(S_\varphi^{1/2}\psi_k\big)_{k\in X}$ in the sense of \cite{Casazza2004}. By \cite[Theorem~2]{Casazza2004} we know that $(\omega_k)_{k\in X}$ is a Riesz sequence if and only if $\big(S_\varphi^{1/2}\psi_k\big)_{k\in X}$ is a frame, which in turn is equivalent to $(\psi_k)_{k\in X}$ being a frame, since $S_\varphi^{1/2}\in \B(\Hil)$ is invertible.\medskip

\textbf{\ref{list:10}$\Rightarrow$\ref{list:9}:} This implication follows directly from Lemma~\ref{lem:invertible}, since we have that $\omega\sim_\A \omega$ by Lemma~\ref{omegalocalization}. 
\medskip

\textbf{\ref{list:9}$ {\Leftrightarrow}$\ref{list:8}:} This follows from duality via Corollary~\ref{SGadoints} (ii).\medskip

\textbf{\ref{list:8}$\Rightarrow$\ref{list:7}:}  
Since $G_\omega = C_{\omega} D_{\omega}:\ell^1(X)\to \ell^1(X)$ is invertible, it follows that $C_\omega:\Hil^1(\varphi)\to \ell^1(X)$ is surjective. 
Furthermore, we have for all $c\in \ell^1(X)$ that 
$$\Vert c \Vert_{\ell^1(X)} = \Vert G_\omega^{-1} C_{\omega} D_{\omega} c \Vert_{\ell^1(X)} \leq \Vert G_\omega^{-1} C_{\omega} \Vert_{\B(\Hil^1(\varphi),\ell^1(X))} \Vert D_{\omega} c \Vert_{\Hil^1(\varphi)}.$$
This means that $D_{\omega}: \ell^1(X) \to \Hil^1(\varphi)$ is bounded from below and therefore has closed range.    
%

Finally, the last claim of Theorem~\ref{thm:main} follows from Remark~\ref{indep}.
\end{proof}

\section{A Comparison to "Gabor Frames Without Inequalities"}\label{A Comparison to "Gabor Frames Without Inequalities"}

As already mentioned in the introduction, the main insight of \cite{gro07ineq} (explicitly stated as a corollary therein) was that, for a window function $g$ in the modulation space $M^1(\R^d)$, the Gabor system $\mathcal{G}=\G(g,\Lambda)$ over a full-rank lattice $\Lambda \subset \mathbb{R}^{2d}$ being a frame for $L^2(\mathbb{R}^d)$ is equivalent to the injectivity of the analysis operator $C_{\G}:M^\infty(\R^d) \to \ell^{\infty}(\Lambda)$. In contrast to Theorem~\ref{thm:main}, no closed range condition is listed in this or any other of the conditions from \cite[Theorem~3.1]{gro07ineq}. This is due to the fact that they are automatically fulfilled in this setting, as the next example illustrates. 

\begin{example}\label{ex:gabor}
For $g\in M^1(\mathbb{R}^d)$ and $\Lambda \subset \mathbb{R}^{2d}$ a full-rank lattice, the following are equivalent:
\begin{itemize}
    \item[(1)] $\G =\G(g,\Lambda)$ is a Gabor frame for $L^2(\mathbb{R}^d)$.
    \item[(2)] $C_{\G}:M^\infty(\R^d) \to \ell^{\infty}(\Lambda)$ is injective and has closed range.
    \item[(3)] $C_{\G}:M^\infty(\R^d) \to \ell^{\infty}(\Lambda)$ is injective.
\end{itemize}
Indeed, (2)$\Rightarrow$(3) is trivial, and (3)$\Rightarrow$(1) was proven in \cite{gro07ineq}. To show (1)$\Rightarrow$(2), we note that due to $g\in M^1(\mathbb{R}^d)$, 
$S_{\G} = D_{\G}C_{\G}$ is bounded and invertible on $M^p(\mathbb{R}^d)$ for all $1\leq p\leq \infty$ \cite[Theorem~3.1]{gro07ineq}. From the latter, we may deduce the desired closed range condition analogously as in the proof of Theorem~\ref{thm:main}.

In a similar manner, other closed range conditions from \cite[Theorem~3.1]{gro07ineq} can be derived.
\end{example}

At the same time, however, all closed range conditions appearing in Theorem~\ref{thm:main} are in fact necessary. To demonstrate this, we "naively" transfer the items (iii)-(vi), (ix), (x) and (xiii) from \cite[Theorem~3.1]{gro07ineq} to our setting, obtaining the following list of conditions:
%
\begin{enumerate} 
    \item[$(a)$] $S_\psi:\H^\infty(\varphi)\to \H^\infty(\varphi)$ is injective,
    \item[$(b)$] $C_\psi:\H^\infty(\varphi)\to \ell^\infty(X)$ is injective,
    \item[$(c)$] $D_\psi:\ell^1(X)\to \H^1(\varphi)$ has dense range,
    \item[$(d)$] $D_\omega:\ell^\infty(X)\to \H^\infty(\varphi)$ is injective,
    \item[$(e)$] $C_\omega:\H^1(\varphi)\to \ell^1(X)$ has dense range,
    \item[$(f)$] $G_\omega:\ell^\infty(X)\to\ell^\infty(X)$ is injective.
\end{enumerate} 
\noindent 
Below, we provide an elementary example of families $\psi$ and $\omega$ in $\Hil$, that are localized with respect to a self-localized Riesz basis $\varphi$, satisfy all the conditions (a)$-$(f), but are not a frame, respectively not a Riesz sequence. 

\begin{example}\label{counterexample}
Let $X=\mathbb{N}$, $\A$ be any uniformly $p$-stable spectral algebra, $\varphi$ be an $\mathcal{A}$-localized Riesz basis (e.g., an orthonormal basis) and define $\psi=(\psi_k)_{k\in \mathbb{N}}:=\big(\frac{1}{k}\widetilde\varphi_k\big)_{k\in \mathbb{N}}$. Then $\omega_k=\frac{1}{k}S_\varphi^{-1/2}\varphi_k$ for each $k\in \mathbb{N}$, and the following hold:
\begin{itemize}
    \item[(i)] One has $\psi \sim_{\A} \varphi$ and, consequently, $\omega \sim_{\A} \varphi$;
    \item[(ii)] $\psi$ is not a frame and, consequently, $\omega$ not a Riesz sequence in $\Hil$;
    \item[(iii)] Items (a)-(f) from above are satisfied;
    \item[(iv)] The range of $C_{\omega}:\Hil^1(\varphi)\to \ell^1(\mathbb{N})$ is not closed in $\ell^1(\mathbb{N})$. 
\end{itemize}

\begin{proof}
$(i)$ It holds $\psi \sim_{\A}\varphi$, because $G_{\psi,\varphi} = \big[\frac{1}{k}\delta_{k,l}\big]_{k,l\in \mathbb{N}}$ is contained in $\A$ due to solidity of $\A$, since $\mathcal{I}_{\A} = [\delta_{k,l}]_{k,l\in \mathbb{N}} \in \A$. The fact that $\omega \sim_{\A} \varphi$ is shown analogously or via Lemma~\ref{omegalocalization}.

$(ii)$ By testing the frame inequalities on $\varphi_l$ one immediately sees that the lower frame inequality is violated, that is, $\psi$ is \emph{not} a frame. Therefore, $\omega$ cannot be a Riesz sequence according to Theorem~\ref{thm:main}. 

$(iii)$ The Gram matrix $G_\omega = \big[\frac{1}{k^2}\delta_{k,l}\big]_{k,l\in \mathbb{N}}:\ell^\infty(\mathbb{N})\to\ell^\infty(\mathbb{N})$ is obviously injective, that is, (f) is satisfied. Consequently, also (d) has to be satisfied, and (e) follows from (d) via duality upon an application of Lemma~\ref{Cadjoint}. 

Next, we show that (a) is satisfied. Let $f=[(f_n)_{n\in \mathbb{N}}]_{\sim_{\varphi}}, g=[(g_n)_{n\in \mathbb{N}}]_{\sim_{\varphi}} \in \Hil^{\infty}(\varphi)$ and assume that $S_{\psi} (f-g)=\sum_{l\in \mathbb{N}}\frac{1}{l^2}\langle f-g,\widetilde{\varphi}_l\rangle\widetilde{\varphi}_l = 0$. Then $0=C_{\varphi}S_{\psi}(f-g) = \big(\frac{1}{k^2}\langle f-g, \widetilde{\varphi}_k\rangle\big)_{k\in \mathbb{N}}$. Thus $0 = \langle f-g, \widetilde{\varphi}_k\rangle = \lim_{n\to\infty}\langle f_n - g_n, \widetilde{\varphi}_k\rangle$ for each $k\in \mathbb{N}$ and in particular $(f_n)_{n\in \mathbb{N}}\sim_{\varphi}(g_n)_{n\in \mathbb{N}}$. This means that $f-g=0\in \Hil^{\infty}(\varphi)$ and (a) is shown. Statement (b) follows immediately from (a), and (c) follows from (b) via duality and Lemma~\ref{Cadjoint}.   

$(iv)$ Since we have already shown that $C_{\omega}:\Hil^1(\varphi)\to \ell^1(\mathbb{N})$ has dense range, it suffices to show that $C_{\omega}:\Hil^1(\varphi)\to \ell^1(\mathbb{N})$ is not surjective. But the latter is clear, since $C_{\omega} = \left( \bigoplus_{k\in \mathbb{N}}\frac{1}{k}\right) C_{\varphi} S_{\varphi}^{-1/2}$, where $C_{\varphi} S_{\varphi}^{-1/2}:\Hil^1(\varphi) \to \ell^1(\mathbb{N})$ is bijective (this follows from Corollary~\ref{cor:D-phi-injective}), hence $\text{Ran}_1(C_\omega) = \lbrace (c_k)_{k\in \mathbb{N}}: c_k = a_k/k, (a_k)_{k\in \mathbb{N}} \in \ell^1(X)\rbrace $ which is not equal to $\ell^1(X)$.
\end{proof}
\end{example}

\begin{rem.}\label{counterexample2}
Via duality (Lemma~\ref{Cadjoint}) and the closed range theorem, we may conclude from Example~\ref{counterexample} (iv), that $D_{\omega}:\ell^{\infty}(\mathbb{N})\to \Hil^{\infty}(\varphi)$ does not have closed range as well. Furthermore, arguing as in the proof of \textbf{\ref{list:5}$\Leftrightarrow$\ref{list:6}}, one may show that $D_{\omega}:\ell^{\infty}(\mathbb{N})\to \Hil^{\infty}(\varphi)$ having closed range is actually equivalent to $C_{\psi}: \Hil^{\infty}(\varphi) \to \ell^{\infty}(\mathbb{N})$ having closed range. Thus $C_{\psi}: \Hil^{\infty}(\varphi) \to \ell^{\infty}(\mathbb{N})$ does not have closed range in the above example. 

Similarly, by switching the roles of $\psi$ and $\omega$ in Example~\ref{counterexample}, one readily sees that it is possible to find $\psi$ and $\omega$ which are not a frame (respectively a Riesz sequence), violate the remaining closed range conditions appearing in Theorem~\ref{thm:main}, while meeting all other conditions appearing in the same result. 

Thus, each of the closed range conditions appearing in Theorem~\ref{thm:main} is necessary for the statement to remain true.
\end{rem.}

\section{Shift-Invariant Spaces}\label{Shift-Invariant Spaces}

In this section, we apply Theorem~\ref{thm:main} to the setting of shift-invariant systems. For a comprehensive overview to shift-invariant spaces, we refer the reader to \cite{aldroubi}.

Let $T_k$ denote the translation operator on $L^p(\mathbb{R}^d)$ ($k\in \mathbb{R}^d$), given by $T_k f(t) = f(t-k)$. For a generator $\varphi\in L^2(\R^d)$, let $V^p(\varphi)$ be the \textit{shift-invariant space} (formally) defined by
$$V^p(\varphi)=\left\{\sum_{k\in\Z^d}c_k T_k\varphi:\ c\in \ell^p(\Z^d)\right\} \qquad (1\leq p\leq \infty).$$
We assume that 
\begin{enumerate}
    \item[(i)] $\varphi $ is continuous;
    \item[(ii)] $|\varphi(x)|\lesssim (1+|x|)^{-s}$ for some $s>d$; and
    \item[(iii)] $(T_k\varphi)_{k\in \Z^d}$ forms a Riesz basis for $V^2(\varphi)$.
\end{enumerate}
Then the assumptions (i)$-$(iii) ensure that $V^p(\varphi)$ is a closed subspace of $L^p(\R^d)$ ($1\leq p\leq \infty$) \cite[Section~4]{groe04}, 
and that $V^2(\varphi)$ is a reproducing kernel Hilbert space with kernel $K_x^\varphi$, i.e.
\begin{equation}\label{eq:RKHS}
f(x)=\langle f,K_x^\varphi\rangle, \qquad f\in V^2(\varphi), x\in \R^d. 
\end{equation}
It was shown in \cite{groe04} that 
$$\vert \langle T_k\varphi, T_l\varphi \rangle \vert \lesssim (1+\vert k-l\vert)^{-s}, \qquad k,l\in \mathbb{Z}^d,$$
i.e., the Riesz basis $(T_k\varphi)_{k\in \Z^d}$ is intrinsically $\J_s(\mathbb{Z}^d)$-localized, where $\J_s(\Z^d)$ denotes the Jaffard algebra from Example~\ref{spectralexamples} (1). Furthermore, from \cite{groe04} we also know that the associated co-orbit spaces $\Hil^p\big((T_k\varphi)_{k\in \mathbb{Z}^d}\big)$ coincide with the shift-invariant spaces $V^p(\varphi)$ ($1\leq p\leq \infty$) 
and 
\begin{equation}\label{polyloc}
\vert \langle K_{x}^\varphi, T_l\varphi \rangle \vert = \vert \varphi(x -l) \vert \lesssim (1+|x-l|)^{-s}, \qquad x\in \mathbb{R}^d,l\in \mathbb{Z}^d.
\end{equation}
Let $X=(x_k)_{k\in\Z^d}\subset\R^d$ be discrete. If $\psi=(\psi_k)_{k\in \Z^d} := (K_{x_k}^\varphi)_{k\in \Z^d}$ forms a frame for $V^2(\varphi)$ with frame bounds $0<A\leq B$, then $X$ is called a \textit{stable set of sampling} for $V^2(\varphi)$, since, due to (\ref{eq:RKHS}), the frame inequality \eqref{eq:frame} takes the form
$$A\|f\|_{L^2(\mathbb{R}^d)}^2\leq \sum_{k\in\Z^d}|f(x_k)|^2\leq  B\|f\|_{L^2(\mathbb{R}^d)}^2,\qquad f\in V^2(\varphi),$$
in this case. Similarly, $X$ is called a \textit{stable set of sampling} for $V^p(\varphi)$ ($1\leq p\leq \infty$), if there exist constants $A_p,B_p >0$ such that 
$$A_p\|f\|_{L^p(\mathbb{R}^d)}\leq \Vert (f(x_k))_{k\in \Z^d} \Vert_{\ell^p(\Z^d)} \leq  B_p \|f\|_{L^p(\mathbb{R}^d)} ,\qquad f\in V^p(\varphi) .$$
Now, assume that $X=(x_k)_{k\in \mathbb{Z}^d} \subset \R^d$ is relatively separated and satisfies the mild perturbation condition
\begin{equation}\label{pert}
\vert x_k -k\vert \leq C,\qquad\text{for some }C>0\text{ and every }k\in\Z^d. 
\end{equation}
Then, since 
$$(1+\vert k-l\vert)^s \leq (1+\vert k-x_k\vert)^s (1+\vert x_k-l\vert)^s, \qquad  k,l\in \mathbb{Z}^d,$$
the estimate (\ref{polyloc}) implies
$$\vert \langle \psi_k, T_l\varphi \rangle \vert \lesssim (1+\vert k-x_k\vert)^{s}(1+\vert k-l\vert)^{-s} \lesssim (1+\vert k-l\vert)^{-s}, \qquad k,l\in \mathbb{Z}^d,$$
which means that $\psi$ and $(T_k\varphi)_{k\in \mathbb{Z}^d}$ are mutually $\J_s(\mathbb{Z}^d)$-localized. Since \(\mathcal J_s(\mathbb Z^d)\) is a uniformly $p$-stable spectral algebra by Theorem \ref{thm:p-stable-examples},
all assumptions of Theorem \ref{thm:main} are met in this setting and all of the conditions \textbf{\ref{list:1}$-$\ref{list:10}} provide necessary and sufficient conditions for such a set $X$ being a stable set of sampling. Instead of repeating all these conditions, we will only reformulate conditions \textbf{\ref{list:1}}, \textbf{\ref{list:4}} and \textbf{\ref{list:8}$-$\ref{list:10}}, since they seem more insightful to us in this case than the remaining ones.

Before doing so, we compute the R-dual $\omega$ of $\psi$. Since the frame operator $S_\varphi$ commutes with translations by $k\in \Z^d$, so does $S_\varphi^{-1/2}$, and we see that  $S_\varphi^{-1/2} T_k\varphi=T_kS_\varphi^{-1/2} \varphi=T_k \gamma$ for $\gamma:=S_\varphi^{-1/2} \varphi$ and $k\in\Z^d$. Consequently,
\begin{equation}\label{eq:omega-SIS}
\omega_k=\sum_{l\in \Z^d}\langle \psi_l , T_k\varphi \rangle S_{\varphi}^{- {1}/{2}}T_l\varphi=\sum_{l\in \Z^d}  \overline{\varphi(x_l-k)}  T_l \gamma,\qquad k\in\Z^d.
\end{equation}
Note, that in case $X=(x_k)_{k\in\Z^d}= \Z^d$, the R-dual $\omega$ of $\psi$ is again a family of shifts of a single generator $\eta$ as
$$\omega_k=T_k\sum_{l\in \Z^d} \overline{ \varphi(l-k)}  T_{l-k} \gamma=T_k \sum_{l\in \Z^d}  \overline{\varphi(l)}  T_{l}\gamma=T_k\eta,\qquad k\in\Z^d.$$
Computing the entries of the associated Gram matrix $G_\omega$, we see that 
\begin{flalign}\label{correlationmatrix}
[G_\omega]_{k,n} &= \langle \omega_n, \omega_k \rangle = \left\langle  \sum_{l\in \Z^d}  \overline{\varphi(x_l-n)} S_{\varphi}^{-1/2} T_l \varphi , \sum_{m\in \Z^d}  \overline{\varphi(x_m-k)} S_{\varphi}^{-1/2} T_m \varphi \right\rangle \notag \\
&= \sum_{l\in \Z^d} \varphi(x_l-k) \overline{\varphi(x_l-n)},
\end{flalign}
since $(S_{\varphi}^{-1/2} T_l \varphi)_{l\in \Z^d}$ is an orthonormal basis. 

\begin{theorem}\label{samplingthm}
Assume that $\varphi \in L^2(\mathbb{R}^d)$ satisfies the conditions (i)$-$(iii), and let $(x_k)_{k\in \mathbb{Z}^d} \subset \mathbb{R}^d$ be relatively separated. Let $G_\omega$ be the autocorrelation matrix with entries given by (\ref{correlationmatrix}) and assume that condition (\ref{pert}) holds. Then the following are equivalent.
\begin{itemize}
\item[(a)] $(x_k)_{k\in \mathbb{Z}^d}$ is a stable set of sampling for $V^2(\varphi)$.
\item[(b)] $(x_k)_{k\in \mathbb{Z}^d}$ is a stable set of sampling for $V^{\infty}(\varphi)$ and $D_{\psi}:\ell^{\infty} \to V^{\infty}(\varphi)$ has closed range. 
\item[(c)] $G_\omega$ is invertible on $\ell^1(\Z^d)$.
\item[(d)] $G_\omega$ is invertible on $\ell^{\infty}(\Z^d)$.
\item[(e)] $G_\omega$ is invertible on $\ell^{2}(\Z^d)$.
\end{itemize} 
\end{theorem}

\begin{rem.}
In the language of \cite{grrost18}, the assumptions (i)$-$(iii) guarantee that $\varphi$ is contained in the Wiener amalgam space $W_0 = W(C,\ell^1)$ and has \emph{stable integer shifts}, which, by \cite[Theorem~3.1]{grrost18}, yields several equivalent conditions for a relatively separated set $(x_k)_{k\in \Z} \subset \mathbb{R}$ being a stable set of sampling for $V^2(\varphi)$. In particular, using a non-commutative version of Wiener's lemma (\cite[Proposition~8.1]{GROCHENIG2015388}), the authors of \cite{grrost18} showed that $(x_k)_{k\in \Z^d}$ being a stable set of sampling for $V^p(\varphi)$ for \emph{some} $1\leq p \leq \infty$ is equivalent to $(x_k)_{k\in \Z^d}$ being a stable set of sampling for $V^p(\varphi)$ for \emph{all} $1\leq p \leq \infty$ (hence the closed range condition of $D_{\psi}:\ell^{\infty} \to V^{\infty}(\varphi)$ from (b) can actually be omitted). Connected to that, we assume condition (\ref{pert}) in order to guarantee localization with respect to the Jaffard algebra $\J_s(\mathbb{Z}^d)$ and utilize inverse-closedness in $\B(\ell^2(\mathbb{Z}^d))$, which can be viewed as another non-commutative version of Wiener's lemma. 

We also note that the autocorrelation matrix $G_\omega$ is given by the product $G_\omega = P_X(\varphi)^* P_X(\varphi)$, where $P_X(\varphi) = (\varphi(x_l -k))_{l,k\in \Z}$ is the matrix appearing in \cite[Theorem~3.1]{grrost18}. According to \cite[Theorem~3.1]{grrost18} (in the case $d=1$), stable sampling is equivalent to $P_\Gamma(\varphi):\ell^{\infty}(\mathbb{Z}^d) \to \ell^{\infty}(\Gamma)$ being bounded and injective for \emph{all} weak limits $\Gamma$ of integer translates of $(x_k)_{k\in \Z}$. In contrast, conditions (c)$-$(e) from above only concern \emph{one} specific operator $G_\omega$. To the best of our knowledge, these invertibility conditions on the autocorrelation matrix $G_\omega$ are new.
\end{rem.}

We conclude this section with a concrete application of Theorem \ref{samplingthm} for $d=1$.

\begin{example}\label{Bsplineex}
Let $\beta(x) = (1-\vert x \vert)_+ = \max(1-\vert x \vert, 0)$ be the centered linear B-spline. Clearly, $\beta$ is continuous and has polynomial decay of any order $s>1$. According to \cite[Theorem~9.2.6]{ole1n}, $(T_k \beta)_{k\in \mathbb{Z}}$ is a Riesz sequence in $L^2(\mathbb{R})$, hence a Riesz basis for $V^2(\beta)$, i.e., $\beta$ is a generator satisfying (i)$-$(iii). Furthermore, for any $0\leq \tau <1$, the relatively separated set $\mathbb{Z} - \tau$ satisfies condition (\ref{pert}) so that all assumptions of Theorem \ref{samplingthm} are met.

\end{example}

\begin{proposition}\label{beta} $\mathbb{Z} - \tau$ is a stable set of sampling for $V^2(\beta)$ if and only if $\tau \neq \frac{1}{2}$.
\end{proposition}

\begin{proof} We check condition (e) of Theorem \ref{samplingthm}. Via (\ref{correlationmatrix}), we see that
$$[G_\omega]_{k,n} = \sum_{l\in \mathbb{Z}} (1-\vert l - \tau -k\vert)_+ \cdot (1-\vert l - \tau -n\vert)_+ = \begin{cases}
(1-\tau)^2 + \tau^2 \quad &(k=n)\\
(1-\tau)\tau \quad &(\vert k-n\vert = 1)\\
0 \qquad &\text{(else)}
\end{cases}$$
due to support reasons. Thus, $G_\omega$ is an infinite Toeplitz matrix, i.e., a convolution operator with respect to the sequence $g = (g_k)_{k\in \mathbb{Z}}$ with entries corresponding to the constant diagonals of $G_\omega$. 
By unitary equivalence via the Fourier transform on $\mathbb{Z}$, $G_\omega$ is invertible on $\ell^2(\mathbb{Z})$ if and only if $\widehat{g}$ is nowhere vanishing. We compute 
\begin{flalign}
\widehat{g}(\xi) &= (1-\tau)\tau e^{2\pi i \xi} + (1-\tau)^2 + \tau^2 + (1-\tau)\tau e^{- 2\pi i \xi} \notag \\
&= (1-\tau)^2 + \tau^2 + 2(1-\tau)\tau \cos(2\pi \xi) \notag \\
&\geq (1-\tau)^2 + \tau^2 - 2(1-\tau)\tau \notag \\
&= (2\tau - 1)^2, \notag
\end{flalign}
and see that $\widehat{g}$ is non-vanishing precisely when $\tau \neq \frac{1}{2}$.
\end{proof} 

\begin{rem.}
In the case $\tau = \frac{1}{2}$ above, $f(x) = \sum_{k\in \mathbb{Z}}(-1)^k (1-\vert x-k \vert)_+ \in V^{\infty}(\beta)$ 
is contained in the kernel of the analysis operator $C_\psi: V^{\infty}(\beta) \to \ell^{\infty}(\mathbb{Z})$ associated to the family $\psi = (K^{\beta}_l)_{l\in \mathbb{Z} - \frac{1}{2}}$, because $f(l-\frac{1}{2}) = (-1)^l \frac{1}{2} + (-1)^{l+1}\frac{1}{2} = 0$
(again due to support reasons), so that $C_{\psi} f = (f(l-\frac{1}{2}))_{l\in \mathbb{Z}} = 0$. So one may also apply Theorem~\ref{thm:main} (4) to get that $\psi$ is not  a frame.

In other words: As in the example from the introduction, a concrete non-trivial kernel element from the space $\Hil^{\infty}(\varphi) = V^\infty(\beta)$ makes the underlying obstruction in Example \ref{Bsplineex} immediately transparent.
\end{rem.}


\section*{Appendix}

We verify Theorem \ref{thm:p-stable-examples}, which says that all spectral algebras listed in
Example~\ref{spectralexamples} are uniformly $p$-stable. Its proof relies on two results from Shin-Sun and Tessera respectively. 

Let \(\mathcal{C}(X)\) be the class of scalar matrices \(A=[A_{k,l}]_{k,l \in X}\) satisfying
\[
 \|A\|_{\mathcal{C}(X)} := \sum_{m\in\mathbb{Z}^d}
 \sup_{\substack{k,l\in X \\
                  k-l\in m+[0,1)^d}}
 |A_{k,l}| <\infty,
\]
where the supremum over the empty set is understood as zero.

The following result is due to Shin and Sun, see \cite[Theorem~2.1]{ShiSun09}

\begin{theorem}\label{thm:shin-sun}\cite{ShiSun09}
Let \(X\subseteq\mathbb R^d\) be relatively separated. Then \(\mathcal{C}(X)\) is uniformly $p$-stable. 
\end{theorem}

In order to show uniform $p$-stability of the Schur algebra, we invoke a result of Tessera, see \cite[Corollary~1.7]{TESSERA20102793}. Following \cite{TESSERA20102793}, a countable discrete metric space $(Y,d_Y)$ is called \emph{doubling with respect to counting measure}, if there exists $D>0$ such that
\begin{equation}\label{doubling}
\# B_{Y}(z,2r) \leq D \# B_Y(z,r)  \qquad (z\in Y, r>0),
\end{equation}
where $B_Y(z,r) = \lbrace y\in Y: d_Y(z,y) \leq r \rbrace$. Now, for $\delta >0$, define the weight $w_\delta(y,y') := 1 + d_Y(y,y')^\delta$ on $Y\times Y$ and let $\mathcal{S}^1_{w_\delta}(Y)$ be the corresponding weighted Schur class, i.e. the class of all scalar matrices $A=[A_{y,y'}]_{y,y'\in Y}$ satisfying
$$\max\left\{ \sup_{y\in Y}\sum_{y'\in Y} \vert A_{y,y'}\vert w_\delta(y,y'),\, \sup_{y'\in Y}\sum_{y\in Y} \vert A_{y,y'}\vert w_\delta(y,y') \right\} <\infty.$$

\begin{theorem}\label{prop:tessera-p-stability}\cite{TESSERA20102793}
Let $(Y,d_Y)$ be a countable discrete metric space whose counting measure is doubling. Then $\mathcal{S}^1_{w_\delta}(Y)$ is uniformly $p$-stable.
\end{theorem}

We cannot directly apply Theorem \ref{prop:tessera-p-stability} to the weighted Schur algebra $\mathcal{S}^1_{\nu}(X)$, since a relatively separated set $X$ need not be doubling with respect to counting measure, see \cite[Remark~2.3]{Sun2008FramesIS}. However, we may circumvent this issue by adding a lattice layer that provides the missing lower volume growth, see below. Let $\vert \cdot \vert_\infty$ denote the maximum-norm on $\mathbb{R}^d$.

\begin{lem.}
\label{lem:doubling-extension}
Let \(X\subseteq\mathbb R^d\) be relatively separated and define
$$Y:= (X\times\{0\})\cup(\mathbb Z^d\times\{1\})$$
with the metric 
$$d_Y((k,i),(l,j)) :=\vert k-l\vert_\infty + \vert i-j \vert.$$
Then $(Y,d_Y)$ is doubling with respect to the counting measure.
\end{lem.}

\begin{proof}
Since both \(X\) and \(\mathbb Z^d\) are relatively separated subsets of $\mathbb{R}^d$, so is $Z:= X\cup\mathbb{Z}^d$. By \cite[Lemma~II.8]{KoBaoff25} applied to $Z$ (with any fixed
\(s_0>d\) and \(\tau_0=\frac{1}{2}\)), there exists \(C>0\) such that
\begin{equation}\label{taulemma}
\# \big(Z\cap B_\infty(z,r) \big)\leq Cr^d, \qquad z\in Z, r \geq 1,    
\end{equation}
where $B_\infty(z,r) = \lbrace z' \in \mathbb{R}^d : \vert z' - z\vert_{\infty} \leq r \rbrace$. Next, let $\iota: Y\to Z,\ \iota(k,i):=k$ be the first coordinate function. Then $\iota$ is surjective and its preimage satisfies
\begin{equation}\label{fiber}
1 \leq \# \iota^{-1}(k) \leq 2, \qquad k\in Z.     
\end{equation}
Now, observe that for \(z'\in B_Y(z,r)\) it holds $\vert \iota(z')-\iota(z)\vert_\infty \leq d_Y(z',z) \leq r$, hence $\iota(B_Y(z,r)) \subseteq Z \cap B_{\infty}(\iota(z),r)$. Together with (\ref{taulemma}) and (\ref{fiber}) this implies
\begin{equation}\label{upper}
\# B_Y(z,r) \leq 2 \# \iota(B_Y(z,r))  \leq  2 \#\big( Z \cap B_\infty(\iota(z),r)\big) \leq 2C r^d, \qquad r\geq 1.
\end{equation}
On the other hand, fix some arbitrary $z = (k,i) \in Y$ and let $r\geq 4$. Choose \(k_0\in\mathbb Z^d\) such that $\vert k-k_0\vert_\infty\leq \frac{1}{2}$ and set $L:=\left\lfloor r-\frac{3}{2}\right\rfloor$. For each \(n\in\mathbb Z^d\) with $\vert n \vert_{\infty} \leq L$  we have 
$$d_Y(z,(k_0 + n,1)) = \vert k - k_0 - n\vert_{\infty} + \vert i -1\vert \leq \frac{1}{2} + L + 1 \leq r.$$
Therefore, each point $(k_0+n,1) \in Y$, where $\vert n \vert_{\infty} \leq L$, is contained in $B_Y(z,r)$. Since there are $(2L+1)^d$-many such points, this implies $\# B_Y(z,r)  \geq (2L+1)^d$. Together with $L = \left\lfloor r-\frac{3}{2}\right\rfloor \geq r - \frac{5}{2}$, we obtain 
$$\# B_Y(z,r)   \geq (2r-4)^d \geq r^d, \qquad r\geq 4.$$
In combination with (\ref{upper}) the latter yields 
\begin{equation}\label{rlarger4}
\# B_Y(z,2r) \leq 2C(2r)^d \leq 2^{d+1}C \# B_Y(z,r) , \qquad z\in Y, r\geq 4,  
\end{equation}
For \(0<r<4\) we see again via (\ref{upper}) that 
$$\# B_Y(z,2r)   \leq \# B_Y(z,8) \leq 2C8^d \leq 2^{3d+1}C \# B_Y(z,r),$$
where we have used $\# B_Y(z,r)  \geq1$. Thus, the doubling estimate (\ref{doubling}) holds for all $z\in Y$ and all \(r>0\) with $D=2^{3d+1}C$.
\end{proof}

\begin{proof}[\textbf{Proof of Theorem~\ref{thm:p-stable-examples}}]
We consider each spectral algebra in
Example~\ref{spectralexamples} separately.

\emph{(1) The Jaffard algebra.}
By Theorem \ref{thm:shin-sun}, it suffices to show that $\mathcal{J}_s(X) \subseteq \mathcal C(X)$ whenever $s>d$. To this end, let \(A = [A_{k,l}]_{k,l\in X} \in\mathcal J_s(X)\). If \(k-l\in m+[0,1)^d\), then $1+\vert m\vert_\infty \leq 2(1+\vert k-l\vert)$. Consequently, for each $k,l\in X$,
$$\vert A_{k,l}\vert \leq \Vert A\Vert_{\mathcal J_s}(1+\vert k-l \vert)^{-s} \leq 2^s\Vert A\Vert_{\mathcal J_s} (1+\vert m\vert_\infty)^{-s}.$$
It follows that
$$ \Vert A\Vert_{\mathcal{C}(X)} \leq 2^s\Vert A\Vert_{\mathcal J_s} \sum_{m\in\mathbb Z^d}(1+\vert m \vert_\infty)^{-s} < \infty ,$$
where the last sum converges due to $s>d$. Hence, $A\in \mathcal{C}(X)$.

\emph{(2) The weighted Schur algebra.}
Let \(A = [A_{k,l}]_{k,l\in X} \in\mathcal S^1_\nu(X)\), and recall in particular the weak growth condition (2c) from Example \ref{spectralexamples}, that is, there exists \(\delta\in(0,1]\) and \(C_\nu>0\) such that
$$\nu(x)\geq C_\nu(1+\vert x \vert)^\delta, \qquad x\in\mathbb R^d.$$
Define \((Y,d_Y)\) as in Lemma~\ref{lem:doubling-extension} and set $w_\delta (z',z) := 1+ d_Y(z',z)^\delta$. We define a matrix \(\widehat A\) indexed by \(Y\times Y\) via 
$$\widehat A_{(k,i),(l,j)} := \begin{cases}
   A_{k,l} \quad &\text{if } i=j=0,\\
   \delta_{k,l} \quad &\text{if } i=j=1,\\
   0 \quad &\text{if } i\neq j.
 \end{cases}$$
Equivalently, we have $\widehat A = A\oplus I_{\ell^p(\mathbb Z^d)}$ under the canonical identification $\ell^p(Y) = \ell^p(X) \oplus \ell^p(\mathbb Z^d)$ ($1\leq p \leq \infty$). For \(k,l\in X\), we estimate
$$w_\delta((k,0),(l,0)) = 1+\vert k-l \vert_\infty^\delta \leq 1+\vert k-l \vert^\delta \leq 2(1+\vert k-l\vert)^\delta \leq \frac{2}{C_{\nu}} \nu(k-l).$$
Hence, we essentially may estimate the $w_\delta$-weighted row- and column sums of $\widehat{A}$ by the corresponding $\nu$-weighted row- and column sums of $A \in \mathcal S^1_\nu(X)$; the rows and columns belonging to the $\mathbb{Z}^d$-block contain only one non-zero entry of weight one and all mixed entries vanish. More precisely, a direct calculation shows 
$$\big\Vert \widehat{A} \big\Vert_{\mathcal S^1_{w_\delta}(Y)} \leq \max \left\lbrace \frac{2}{C_{\nu}} \Vert A\Vert_{\mathcal S^1_\nu(X)} \, ,  1 \right\rbrace < \infty$$
so that $\widehat A\in \mathcal S^1_{w_\delta}(Y)$. Next, we claim that for every \(p\in[1,\infty]\),
\begin{equation}\label{lambdapestimate}
\lambda_p \big( \widehat A \, \big) = \min \{\lambda_p(A),1\}.    
\end{equation}
Indeed, for $\widehat{c} = (c,d) \in \ell^p(Y)$, where $c\in \ell^p(X), d\in \ell^p(\mathbb{Z}^d)$, we have by definition of $\widehat{A}$ that $\widehat{A} \widehat{c} = \widehat{A}(c,d) = (Ac,d) \in \ell^p(Y)$. Set \(m_p:=\min\{\lambda_p(A),1\}\). If
\(1\leq p<\infty\), then
\[
\begin{aligned}
 \Vert\widehat A \widehat{c} \Vert_{\ell^p(Y)}^p  &=\Vert Ac\Vert_{\ell^p(X)}^p + \Vert d\Vert_{\ell^p(\mathbb{Z}^d)}^p \\
 &\geq \lambda_p(A)^p \Vert c\Vert_{\ell^p(X)}^p + \Vert d\Vert_{\ell^p(\mathbb{Z}^d)}^p \\
 &\geq m_p^p \bigl(\Vert c\Vert_{\ell^p(X)}^p + \Vert d\Vert_{\ell^p(\mathbb{Z}^d)}^p\bigr) = m_p^p \Vert \widehat{c} \Vert_{\ell^p(Y)}^p.
\end{aligned}
\]
Thus $\lambda_p(\widehat A) \geq \min \{\lambda_p(A),1\}$. 
Conversely, testing \(\widehat A\) on vectors of the form
\((c,0)\) and \((0,d)\), respectively, gives $\lambda_p (\widehat A) \leq \lambda_p(A)$ and $\lambda_p(\widehat A)\leq1$, respectively. In the case $p=\infty$ one argues similarly, whence the claim follows. Finally, suppose now that \(\lambda_p(A)>0\) for some \(1\leq p\leq \infty\). By (\ref{lambdapestimate}), \(\lambda_p(\widehat A)>0\). Since, by Lemma \ref{lem:doubling-extension}, $(Y,d_Y)$ is doubling with respect to counting measure and since $\widehat{A} \in S^1_{w_\delta}(Y)$, Theorem~\ref{prop:tessera-p-stability} gives \(\lambda_q(\widehat A)>0\) for all \( 1 \leq q \leq \infty\). Another application of (\ref{lambdapestimate}) yields \(\lambda_q(A)>0\) for every \(1\leq q \leq \infty\). Thus \(\mathcal S^1_\nu(X)\) is uniformly $p$-stable.

\emph{(3) The weighted Baskakov-Gohberg-Sj\"ostrand algebra.}
Let \(A\in\mathcal C_\nu(\mathbb Z^d)\), where $\nu$ is submultiplicative, symmetric and satisfies the GRS-condition. Then we have in particular that $\nu \geq 1$ on $\mathbb{R}^d$ so that 
$\Vert A\Vert_{\mathcal{C}(\mathbb{Z}^d)} \leq \Vert A\Vert_{\mathcal{C}_{\nu}(\mathbb{Z}^d)}$. Thus, $     \mathcal C_\nu(\mathbb Z^d) \subseteq \mathcal C(\mathbb Z^d)$ and uniform $p$-stability follows again from Theorem~\ref{thm:shin-sun} applied to $X=\mathbb{Z}^d$.

\emph{(4) The anisotropic algebras.}
The anisotropic versions of the algebras considered in (1)--(3), see \cite[Section~6]{KoeBa25}, are subalgebras of their respective isotropic counterparts. Since uniform $p$-stability is inherited by subalgebras, their uniform $p$-stability follows from the uniform $p$-stability established in
(1)--(3).
\end{proof}

\section*{Acknowledgments}

\noindent The authors thank K. Gr\"ochenig for pointing out the connection between Theorem \ref{samplingthm} and \cite{grrost18}, and also  O. Christensen and E. King for their helpful comments.

The authors thank the anonymous reviewer for valuable comments that helped to improve the quality of this article.


The authors used ChatGPT 5.6 Sol (OpenAI) and DeepL for language translation,  editorial assistance and literature search. The idea for the concrete B-spline application of Theorem~\ref{samplingthm} appearing in Example \ref{Bsplineex} and Proposition \ref{beta}, as well as the proof of the uniform $p$-stability of $\MS_\nu^1$ via Lemma \ref{lem:doubling-extension} was suggested by ChatGPT and subsequently independently summarized and verified by all authors. All remaining mathematical content was independently developed and verified by the authors.

This research was funded by the Austrian Science Fund (FWF) 10.55776/P34624 (P.B. and L.K.) and 10.55776/PAT1384824 (M.S.).
For open
access purposes, the authors have applied a CC BY public copyright license to any author-accepted manuscript
version arising from this submission.

\bibliographystyle{abbrv} 
\bibliography{biblioall}

@article{groelei06,
 ISSN = {00029947},
 URL = {http://www.jstor.org/stable/3845503},
 author = {Karlheinz Gr{\"o}chenig and Michael Leinert},
 journal = {Trans. Amer. Math. Soc.},
 number = {6},
 pages = {2695--2711},
 publisher = {American Mathematical Society},
 title = {Symmetry and Inverse-Closedness of Matrix Algebras and Functional Calculus for Infinite Matrices},
 urldate = {2023-09-08},
 volume = {358},
 year = {2006}
}

@article{barnes87,
 ISSN = {03794024, 18417744},
 URL = {http://www.jstor.org/stable/24714411},
 author = {Bruce A. Barnes},
 journal = {J. Oper. Theory},
 number = {1},
 pages = {115--132},
 publisher = {Theta Foundation},
 title = {The spectrum of integral operators on {L}ebesgue spaces},
 urldate = {2023-09-11},
 volume = {18},
 year = {1987}
}

@Article{all,
  Title                    = {Short term spectral analysis, synthesis, and modification by discrete {F}ourier transform},
  Author                   = {J.B. Allen},
  Journal                  = {{IEEE} Trans. Acoust., Speech, Signal Proces.},
  Year                     = {1977},
  Number                   = {3},
  Pages                    = {235--238},
  Volume                   = {ASSP-25}
}

@inproceedings{KoBaoff25,
  author    = {K{\"o}hldorfer, Lukas and Balazs, Peter},
  title     = {On the Inverse-Closedness of Operator-Valued Matrices with Polynomial Off-Diagonal Decay},
  booktitle = {2025 International Conference on Sampling Theory and Applications (SampTA)},
  pages     = {1--5},
  publisher = {IEEE},
  year      = {2025},
  doi       = {10.1109/SampTA64769.2025.11133549}
}

@INPROCEEDINGS{haubakoe25,
author={Hauschka, Nikolas and Balazs, Peter and K\"ohldorfer, Lukas},
booktitle={2025 International Conference on Sampling Theory and Applications (SampTA)}, 
title={{B}anach distribution spaces for a {H}ilbert space}, 
year={2025},
volume={},
number={},
pages={1-5},
doi={10.1109/SampTA64769.2025.11133572}
}

@book{GRS64,
  title={Commutative {N}ormed {R}ings},
  author={Gelfand, I.M. and Raikov, D.A. and Shilov, G.E.},
  publisher={Chelsea Publishing Company},
  Address={Bronx, New York},
  year={1964},
  url={https://books.google.at/books?id=clIWcgAACAAJ},
}

@InCollection{xxlgro14,
  author    = {Balazs, Peter and Gr{\"o}chenig, Karlheinz},
  title     = {A Guide to Localized Frames and Applications to {G}alerkin-like Representations of Operators},
  booktitle = {Frames and Other Bases in Abstract and Function Spaces},
  year      = {2017},
  editor    = {Isaac Pesenson and Hrushikesh Mhaskar and Azita Mayeli and Quoc T. Le Gia and Ding-Xuan Zhou},
  series    = {Applied and Numerical Harmonic Analysis series (ANHA)},
  publisher = {Birkhauser/Springer},
  url       = {https://arxiv.org/abs/1611.09692},
}

@Article{xxlstoeant11,
  author       = {Peter Balazs and Diana Stoeva and Jean-Pierre Antoine},
  title        = {Classification of {G}eneral {S}equences by {F}rame-{R}elated {O}perators},
  journal      = {Sampl. Theory Signal Image Process.},
  journaltitle = {Sampling Theory in Signal Processing},
  year         = {2011},
  volume       = {10},
  number       = {2},
  pages        = {151-170},
}

@Article{Casazza2004,
  Title                    = {Duality Principles in Frame Theory},
  Author                   = {Casazza, Peter G. and Kutyniok, Gitta and Lammers, Mark C.},
  Journal                  = {J. Fourier Anal. Appl.},
  Year                     = {2004},
  Number                   = {4},
  Pages                    = {383-408},
  Volume                   = {10},

  Doi                      = {10.1007/s00041-004-3024-7},
  ISSN                     = {1069-5869},
  Language                 = {English},
  Publisher                = {Birkh{\"a}user-Verlag},
  Url                      = {http://dx.doi.org/10.1007/s00041-004-3024-7}
}

@Book{ole1n,
  Title                    = {{A}n {I}ntroduction to {F}rames and {R}iesz {B}ases},
  Author                   = {Christensen, Ole},
  Publisher                = {Birkh{\"a}user},
  Year                     = {2016},

  Notes                    = {% Book in personal Library!}
}

@Book{conw1,
  Title                    = {A Course in Functional Analysis},
  Author                   = {J. B. Conway},
  Publisher                = {Springer},
  Year                     = {1990},
  Edition                  = {2.},
  Series                   = {Graduate Texts in Mathematics},
  address                  = {New York},
  Notes                    = {Personal library. Full name: John B. Conway}
}

@article{KoeBa25,
title = {Wiener pairs of {B}anach algebras of operator-valued matrices},
journal = {J. Math. Anal. Appl.},
volume = {549},
number = {2},
pages = {129525},
year = {2025},
issn = {0022-247X},
doi = {https://doi.org/10.1016/j.jmaa.2025.129525},
url = {https://www.sciencedirect.com/science/article/pii/S0022247X25003063},
author = {Lukas K\"ohldorfer and Peter Balazs},
}

@Article{groe04,
  author    = {Karlheinz Gr\"{o}chenig},
  title     = {Localization of frames, {B}anach frames, and the invertibility of the frame operator},
  doi = {https://doi.org/10.1007/s00041-004-8007-1},
  volume    = {10},
  journal   = {J. Fourier Anal. Appl.},
  pages     = {105-132},
  year      = {2004},
}

@Article{grrost18,
  author       = {{G}r\"ochenig, {K}arlheinz and {R}omero, {J}ose-{L}uis and {S}t\"ockler, {J}oachim},
  title        = {Sampling theorems for shift-invariant spaces, {G}abor frames, and totally positive functions},
  pages        = {1119-1148},
  volume       = {211},
  journal      = {{I}nvent. {M}ath.},
  year         = {2018},
}

@Article{forngroech1,
  author  = {Massimo Fornasier and Karlheinz Gr{\"o}chenig},
  title   = {Intrinsic Localization of Frames},
  journal = {Constr. Approx.},
  year    = {2005},
  volume  = {22},
  number  = {3},
  pages   = {395--415},
}

@InBook{gr10-2,
  author       = {{G}r{\"o}chenig, {K}arlheinz},
  title        = {{W}iener's lemma: {T}heme and variations. {A}n introduction to spectral invariance and its applications.},
  booktitle    = {Four Short Courses on Harmonic Analysis. Applied and Numerical Harmonic Analysis},
  year         = {2010}, 
  series       = {{A}pplied and {N}umerical {H}armonic {A}nalysis},
  publisher    = {{B}irkh{\"a}user},
  chapter      = {5},
  pages        = {175 - 234},
  organization = {{N}u{H}{A}{G}},
}

@Article{gr91,
  author       = {{G}r{\"o}chenig, {K}arlheinz},
  title        = {{D}escribing functions: atomic decompositions versus frames},
  number       = {3},
  pages        = {1--41},
  volume       = {112},
  journal      = {{M}onatsh. {M}ath.},
  organization = {{N}u{H}{A}{G};{C}lassical},
  year         = {1991},
}

@Article{grorro15,
  author       = {{G}r{\"o}chenig, {K}arlheinz and {O}rtega-{C}erd{\`a}, {J}oaquim and {R}omero, {J}os{\'e} {L}uis},
  title        = {{D}eformation of {G}abor systems},
  number       = {4},
  pages        = {388-425},
  volume       = {277},
  journal      = {{A}dv. {M}ath.},
  organization = {{L}arfex;locetisa;{N}u{H}{A}{G}},
  year         = {2015},
}

@Article{gro07ineq,
  Title                    = {{G}abor frames without inequalities},
  Author                   = {Karlheinz Gr{\"o}chenig},
  Journal                  = {Int. Math. Res. Not. IMRN},
  Year                     = {2007},
  Number                   = {23},
  Pages                    = {rnm111},
  Volume                   = {2007},

  Link                     = {http://arxiv.org/abs/math.FA/0703379},
  Url                      = {http://arxiv.org/abs/math/0703379}
}

@Book{gr01,
  author    = {Karlheinz Gr{\"o}chenig},
  title     = {Foundations of Time-Frequency Analysis},
  year      = {2001},
  publisher = {Birkh{\"a}user},
  address    = {Boston},
  notes     = {% full name: Karlheinz Gr{\"o}chenig},
}

@article{Sjostrand19941995,
author = {Sj\"ostrand, J.},
journal = {S\'eminaire {\'E}quations aux d{\'e}riv{\'e}es partielles (Polytechnique)},
language = {eng},
pages = {1-19},
publisher = {Ecole Polytechnique, Centre de Mathématiques},
title = {Wiener type algebras of pseudodifferential operators},
url = {http://eudml.org/doc/112115},
year = {1994-1995},
}

@article{aldroubi,
author={Aldroubi, A. and Gr{\"o}chenig, K.},
title={Nonuniform sampling and reconstruction in shift-invariant spaces},
journal={SIAM Rev.},
volume=43,
number=4,
year=2001,
pages={585-620},
}

@Book{Maddox:101881,
  author    = {Ivor J Maddox},
  title     = {{Infinite {M}atrices of Operators}},
  year      = {1980},
  series    = {Lecture Notes in Mathematics},
  publisher = {Springer},
  url       = {https://cds.cern.ch/record/101881},
  address   = {Berlin},
}

@Article{ja90,
  author   = {St{\'e}phane Jaffard},
  title    = {Propri{\'e}t{\'e}s des matrices "bien localis{\'e}es" pre{\'e} de leur diagonale et qualques applications},
  journal  = {Ann. Inst. H. Poincar{\'e} Anal. Non Lin{\'e}aire},
  year     = {1990},
  volume   = {7},
  number   = {5},
  pages    = {461--476},
 }

@article{ron-shen,
    author = {Ron, A. and Shen, Z.},
    title = {Weyl-{H}eisenberg frames and {R}iesz bases in ${L}^2(\mathbb{R}^d)$},
    journal = {Duke Math. J.},
    year = 1997,
volume=89,
number=2,
pages={237-282},
}

@article{bacahela06a,
author = {Balan, Radu and Casazza, Peter and Heil, Christopher and Landau, Zeph},
year = {2006},
month = {04},
pages = {105-143},
title = {Density, {O}vercompleteness, and {L}ocalization of {F}rames. {I}. {T}heory},
volume = {12},
journal = {J. Fourier Anal.  Appl.},
doi = {10.1007/s00041-006-6022-0},
}

@article{dahlkeetal08,
author = {Dahlke, Stephan and Fornasier, Massimo and Rauhut, Holger and Steidl, Gabriele and Teschke, Gerd},
title = {Generalized coorbit theory, {B}anach frames, and the relation to {$\alpha$}-modulation spaces},
journal = {Proc. London Math. Soc.},
volume = {96},
number = {2},
pages = {464--506},
year = {2008},
doi = {10.1112/plms/pdm051}
}

@Article{Balan2006,
  author   = {Balan, Radu and Casazza, Peter G. and Heil, Christopher and Landau, Zeph},
  title    = {Density, Overcompleteness, and Localization of Frames. {I}{I}. {G}abor Systems},
  journal  = {J. Fourier Anal. Appl.},
  year     = {2006},
  volume   = {12},
  number   = {3},
  pages    = {307--344},
  issn     = {1531-5851},
  doi      = {10.1007/s00041-005-5035-4},
  url      = {https://doi.org/10.1007/s00041-005-5035-4},
  day      = {01},
}

@Article{Baskakov1990,
  author  = {Baskakov, A. G.},
  title   = {Wiener's theorem and the asymptotic estimates of the elements of inverse matrices},
  journal = {Funct. Anal. Appl.},
  year    = {1990},
  volume  = {24},
  number  = {3},
  pages   = {222--224},
  issn    = {1573-8485},
  doi     = {10.1007/BF01077964},
  url     = {https://doi.org/10.1007/BF01077964},
  day     = {01},
}

@article{Gohberg89,
author = {Gohberg, I. and Kaashoek, M. and Woerdeman, Hugo},
year = {1989},
month = {01},
pages = {343-382},
title = {The band method for positive and strictly contractive extension problems: An alternative version and new applications},
volume = {12},
journal = {Integral Equ. Oper. Theory},
doi = {10.1007/BF01235737}
}

@Article{xxlgrospeck19,
  author       = {Peter Balazs and Karlheinz Gr{\"o}chenig and Michael Speckbacher},
  journaltitle = {Transactions of the American Mathematical Society Series B},
  title        = {Kernel theorems in coorbit theory},
  doi          = {https://doi.org/10.1090/btran/42},

  pages        = {346-364},
  volume       = {6},
  journal      = {Trans. Am. Math. Soc. Ser. B},
  year         = {2019},
}

@article{sun,
    author = {Sun, Q.},
    title = {Wiener's lemma for infinite matrices},
volume =359, 
number =7, 
year=2007, 
pages ={3099-3123},
    journal = {Trans. Amer. Math. Soc.},
}

@article{xxlbysp24,
    author = {Bytchenkoff, D. and Speckbacher, M. and Balazs, P.},
title={Kernel theorems for operators on co-orbit spaces associated with localised frames},
    journal = {J. Math. Anal. Appl.},
    year = 2025,
volume=551,
number=1,
}

@Article{feichtinger80cras,
  author   = {Feichtinger, H. G.},
  title    = {Un espace de {B}anach de distributions temp\'er\'ees sur les groupes localement compacts ab\'eliens},
  journal  = {C. R. Acad. Sci. Paris S\'er. A-B},
  year     = {1980},
  volume   = {290},
  number   = {17},
  pages    = {A791--A794},
  issn     = {0151-0509},
  coden    = {CHASAP},
  mrclass  = {43A15},
  mrnumber = {81e:43006},
}

@article{GROCHENIG2015388,
title = {Deformation of {G}abor systems},
journal = {Adv. Math.},
volume = {277},
pages = {388-425},
year = {2015},
issn = {0001-8708},
doi = {https://doi.org/10.1016/j.aim.2015.01.019},
url = {https://www.sciencedirect.com/science/article/pii/S0001870815000730},
author = {Karlheinz Gr\"ochenig and Joaquim Ortega-Cerd\`a and José Luis Romero}
}

@Article{B,
  author    = {B. Bodmann},
  title     = {Optimal linear transmission by loss-insensitive packet encoding},
  journal   = {Appl. Comp. Har. Anal.},
  year      = {2007},
  volume    = {22},
  pages     = {274--285},
}

@article{Sun2008FramesIS,
  title={Frames in spaces with finite rate of innovation},
  author={Qiyu Sun},
  journal={Advances in Computational Mathematics},
  year={2008},
  volume={28},
  pages={301-329},
  url={https://api.semanticscholar.org/CorpusID:18119591}
}

@article{ShiSun09,
author = {Shin, Chang and Sun, Qiyu},
year = {2009},
month = {11},
pages = {2417--2439},
title = {Stability of {L}ocalized {O}perators},
volume = {256},
journal = {Journal of Functional Analysis},
doi = {10.1016/j.jfa.2010.07.014}
}

@article{lifting,
    author = {Balazs, P. and Gr\"ochenig, K.},
    title = {The lifting property for frame multipliers and {T}oeplitz operators},
    journal =  	{arXiv:2506.18441},
    year = 2025,
}

@Article{ALDROUBI20081667,
  author   = {Akram Aldroubi and Anatoly Baskakov and Ilya Krishtal},
  title    = {Slanted matrices, {B}anach frames, and sampling},
  journal  = {J. Funct. Anal.},
  year     = {2008},
  volume   = {255},
  number   = {7},
  pages    = {1667 - 1691},
  issn     = {0022-1236},
  doi      = {https://doi.org/10.1016/j.jfa.2008.06.024},
  url      = {http://www.sciencedirect.com/science/article/pii/S0022123608002437},
}

@Article{GrLe04,
  author    = {Karlheinz Gr\"ochenig and Michael Leinert},
  title     = {Wiener's Lemma for Twisted Convolution and {G}abor Frames},
  issn      = {08940347, 10886834},
  number    = {1},
  pages     = {1--18},
  url       = {http://www.jstor.org/stable/20161184},
  volume    = {17},
  journal   = {J. Amer. Math. Soc.},
  year      = {2004},
}

@MastersThesis{hol22,
  author      = {Jakob Holb\"ock},
  school = {University of Vienna},
  title       = {Localized frames and applications},
  year        = {2022},
}

@article{TESSERA20102793,
title = "Left inverses of matrices with polynomial decay",
journal = "Journal of Functional Analysis",
volume = "259",
number = "11",
pages = "2793 - 2813",
year = "2010",
issn = "0022-1236",
doi = "https://doi.org/10.1016/j.jfa.2010.07.014",
url = "http://www.sciencedirect.com/science/article/pii/S0022123610002892",
author = "Romain Tessera"
}

@article{r-duals,
    title = {On {R}-duals and the duality principle in {G}abor analysis},
author={Stoeva, D. T. and Christensen, O.},
journal={J. Fourier Anal. Appl.},
volume=21, 
pages={383-400},
year=2015,
}

@inbook{feilucor08,
author = {Feichtinger, Hans and Luef, Franz and Cordero, E.},
year = {2008},
month = {08},
pages = {1-33},
title = {{B}anach {G}elfand {T}riples for {G}abor {A}nalysis},
volume = {1949},
isbn = {978-3-540-68266-0},
booktitle = {Pseudo-Differential Operators. Lecture Notes in Mathematics},
doi = {10.1007/978-3-540-68268-4_1},
publisher={Springer, Berlin, Heidelberg}
}

@article{koeba24,
    author = "L. K{\"o}hldorfer and P. Balazs",
    title   = "Localization of operator-valued frames",
    year    = "2026",
    journal = "Appl. Comput. Harmon. Anal.",
    volume=84,
    pages=101868,
}

\end{document}